\documentclass[11pt]{amsart}
\usepackage{macros}

\title[Asymptotic normality of the major index on standard tableaux]{Asymptotic normality of the \\ major index on standard tableaux}

\author{Sara C. Billey, Matja\v{z} Konvalinka, Joshua P. Swanson}
\address{Billey: Department of Mathematics, University of Washington,
  Seattle, WA 98195, USA}
\email{billey@math.washington.edu}

\thanks{The first author was partially supported by the Washington
  Research Foundation and  DMS-1764012 from the National Science Foundation. The second author was partially supported by Research Project
BI-US/16-17-042 of the Slovenian Research Agency and research core funding No. P1-0294.}
\address{Konvalinka: Faculty of Mathematics and Physics,
University of Ljubljana,
Jadranska 21, Ljubljana, Slovenia, and Institute for Mathematics, Physics and Mechanics, Jadranska 19, Ljubljana, Slovenia}
\email{matjaz.konvalinka@fmf.uni-lj.si}
\address{Swanson: Department of Mathematics,
University of California, San Diego (UCSD),
La Jolla, CA  92093-0112}
\email{jswanson@ucsd.edu}
\date{\today}

\begin{document}

\begin{abstract}
  We consider the distribution of the major index on
  standard tableaux of arbitrary straight shape and certain skew
  shapes. We use cumulants to classify all possible limit laws
  for any sequence of such shapes in terms of a simple auxiliary
  statistic, $\aft$, generalizing earlier results of
  Canfield--Janson--Zeilberger, Chen--Wang--Wang, and others.
  These results can be interpreted as giving a very precise description of
  the distribution of irreducible representations in different degrees
  of coinvariant algebras of certain complex reflection groups.
  We conclude with some conjectures concerning unimodality,
  log-concavity, and local limit theorems.
\end{abstract}
\maketitle

\setcounter{tocdepth}{1}
\tableofcontents

\section{Introduction}
\label{sec:intro}

The study of permutation and partition statistics is a classic topic
in enumerative combinatorics. The major index
statistic on permutations was introduced a century ago by Percy MacMahon in his
seminal works \cite{MacMahon.1913,MR1576566}.  This statistic, denoted
$\maj(w)$, is defined to be the sum of the positions of the descents
of the permutation $w=[w_1,w_2,\ldots,w_n]$ in one-line notation.  A
descent is any position $i$ such that $w_i>w_{i+1}$.  At first glance,
this function on permutations may be unintuitive, but it has inspired
hundreds of papers and many generalizations; for example on Macdonald
polynomials \cite{HHL2005}, posets \cite{Ehrenborg.Readdy.2015},
quasisymmetric functions \cite{Shareshian-Wachs.2010}, cyclic sieving
\cite{Reiner-Stanton-White.CSP,AHLBACH201837}, and bijective
combinatorics \cite{Foata,carlitz.1975}.

The following central limit theorem for
$\maj$ on $S_n$ is well known and is an archetype for our results.
Given a real-valued random variable $\cX$, we let
\[ \cX^* \coloneqq \frac{\cX - \mu}{\sigma} \]
denote the corresponding normalized random variable with mean $0$
and variance $1$.  Briefly, we say $\maj$ on $S_n$ is
\textit{asymptotically normal} as $n \to \infty$ based on the following
classical result. See \Cref{tab:an} for further examples.

\begin{Theorem}
  \cite{MR0013252}
  \label{thm:maj_Sn_an}
  Let $\cX_n[\maj]$ denote the major index random variable on $S_n$ under the
  uniform distribution. Then, for all $t \in \bR$,
    \[ \lim_{n \to \infty} \bP[\cX_n[\maj]^* \leq t] = \bP[\cN \leq t] \]
  where $\cN$ is the standard normal random variable.
\end{Theorem}

\begin{table}
  \begin{center}
  \begin{tabular}{p{3cm}|p{2.6cm}|p{2.6cm}|p{3.4cm}}
  \toprule
  Statistic & Set & Generating Function & References \\
  \midrule\midrule
      $\#$ elements
      & subsets
      & $(1+q)^n$
      & classical \\
    \midrule
          $\#$ parts
      & strict partitions
      & $\prod_{m=1}^\infty (1 + xy^m)$
      & \cite{MR0004841} \\
    \midrule
          length/inversion number/major index
      & $S_n$
      & $[n]_q!$
      & \cite{MR0013252}, \cite{gon44} \\
    \midrule     $\#$ cycles; $\#$ left-to-right minima
      & $S_n$
      & $\prod_{i=0}^{n-1} (q + i)$
      & \cite{MR0013252}, \cite{gon44} \\
    \midrule
      $\#$ descents
      & $S_n$
      & Eulerian \qquad polynomial $A_n(q)$
      & \cite[pp. 150--154]{MR0155371} \\
    \midrule
      $\#$ descents
      & conjugacy classes in $S_n$
      & \cite[Thm.~1]{MR1652841}
      & \cite{MR1652841,1803.10457} \\
    \midrule
    $\#$ blocks
      & set partitions
      & $\sum_k S(n, k) q^k$
      & \cite{MR0211432} \\
    \midrule
    $\#$ valleys
      & Dyck paths
      & $\frac{1}{[n+1]_q}\binom{2n}{n}_q$
      & \cite[Cor.~3.3]{Chen-Wang-Wang.2008}; \cite[p.~255]{MR814413}
      \\
    \midrule
    length/inversion number/major index
      & $S_n/S_J$, words type $\alpha$
      & $\binom{n}{\alpha}_q$
      & see \Cref{rem:cjz_priority} \\
   \midrule
    major index
      & $\SYT(\lambda)$
      & $q^{\rank(\lambda)} \frac{[n]_q!}{\prod_{c \in \lambda} [h_c]_q}$
      & \Cref{thm:an} \\
   \bottomrule
  \end{tabular}
  \end{center}
  \caption{Summary of some asymptotic normality results for combinatorial
  statistics. See \cite[Ch.~3]{MR3408702}.}
  \label{tab:an}
\end{table}

In this paper, we study the distribution of the major index statistic
generalized to standard Young tableaux of straight and skew shapes.
The properties we discuss here naturally generalize known properties
of the major index distribution on permutations.  They also have
representation theoretic consequences in terms of coinvariant algebras
of complex reflection groups. We will briefly introduce the main results.
See \Cref{sec:background} for more details on the background.

\begin{figure}[ht]
  \centering
  \begin{subfigure}[t]{0.32\textwidth}
    \centering
    \includegraphics[width=\textwidth]{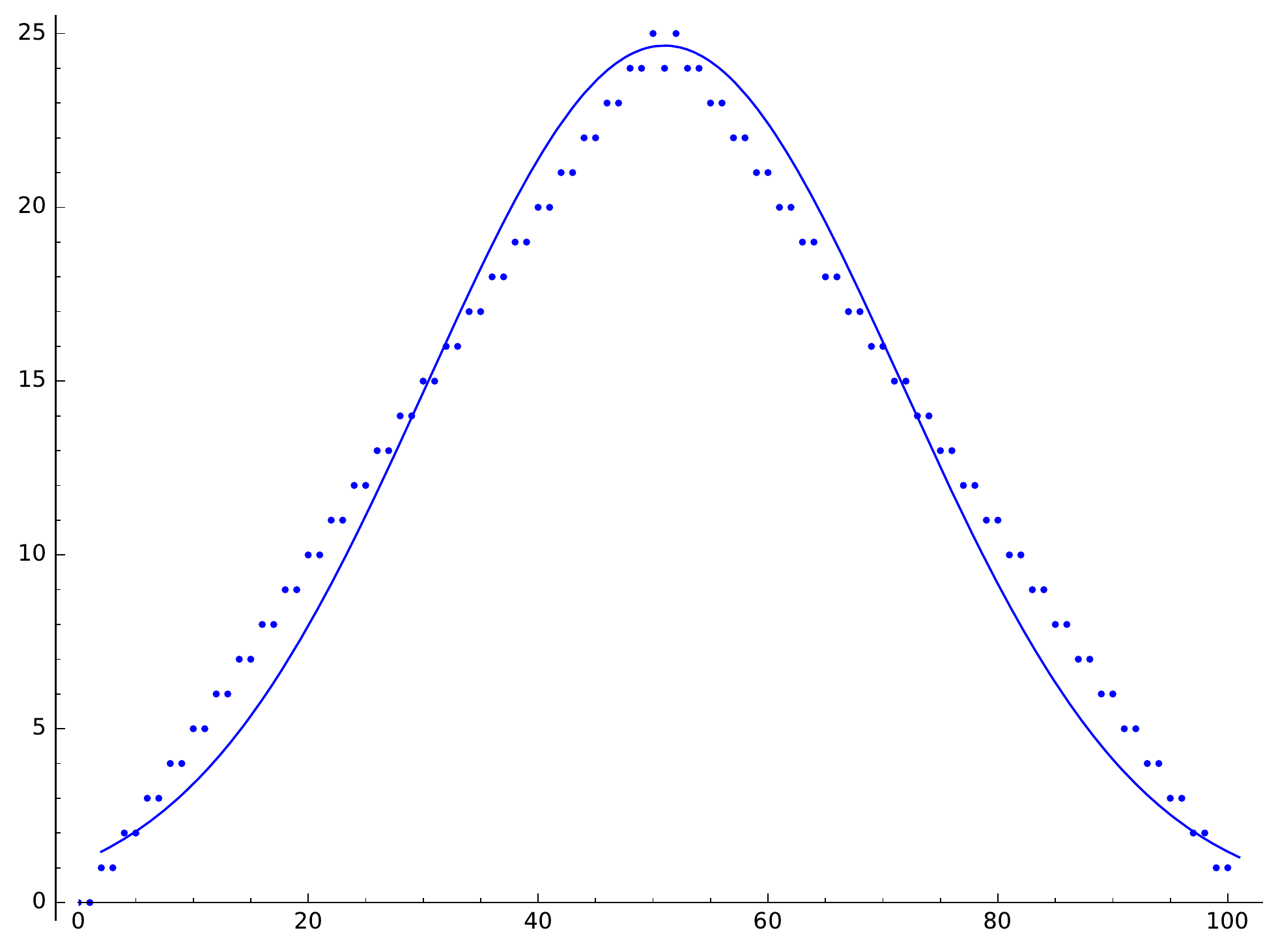}
    \caption{$\lambda = (50, 2)$, $\aft(\lambda) = 2$}
    \label{fig:an_examples_a}
  \end{subfigure}
  \begin{subfigure}[t]{0.32\textwidth}
    \centering
    \includegraphics[width=\textwidth]{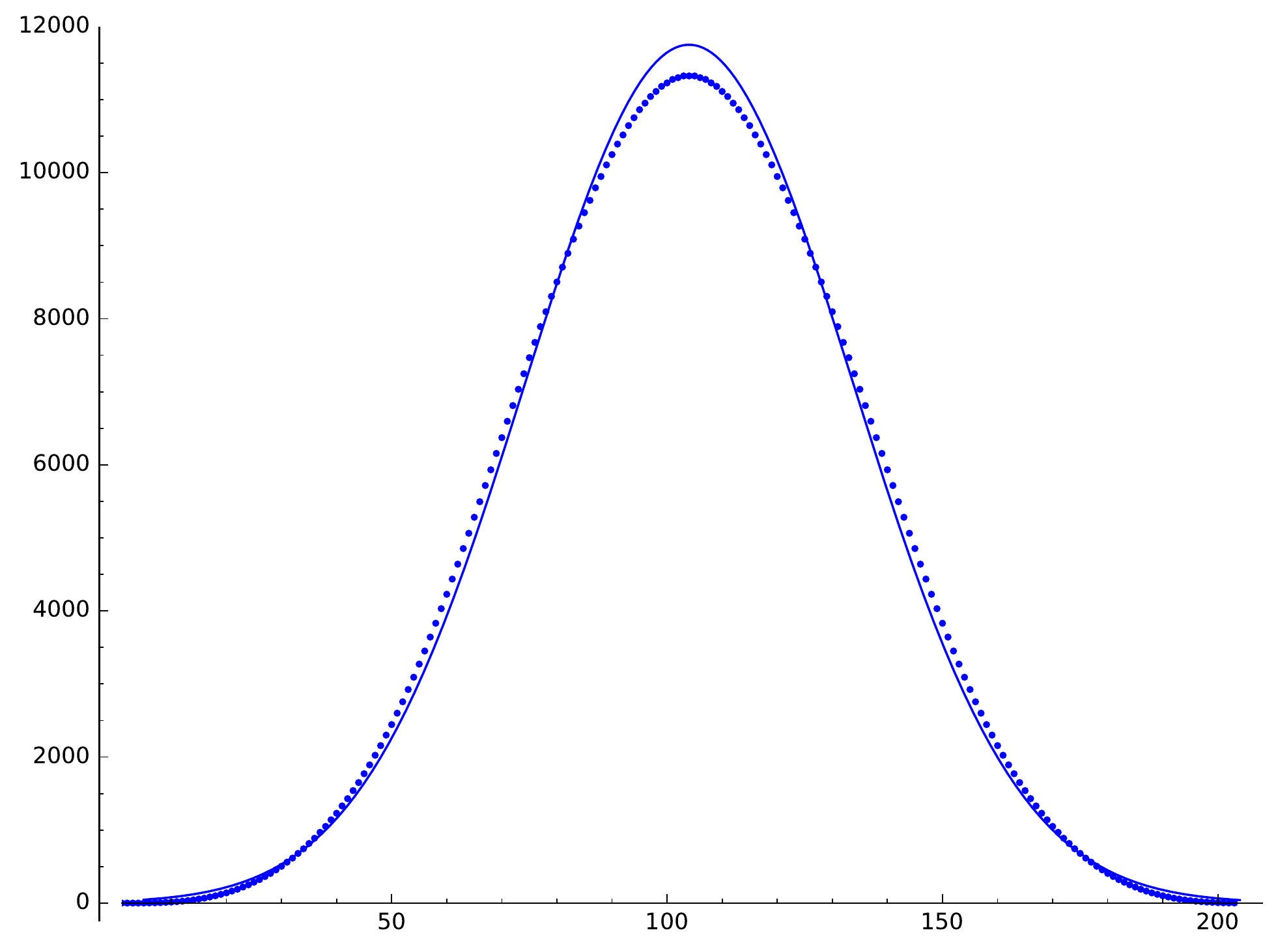}
    \caption{$\lambda = (50, 3, 1)$, $\aft(\lambda) = 4$}
    \label{fig:an_examples_b}
  \end{subfigure}
  \begin{subfigure}[t]{0.32\textwidth}
    \centering
    \includegraphics[width=\textwidth]{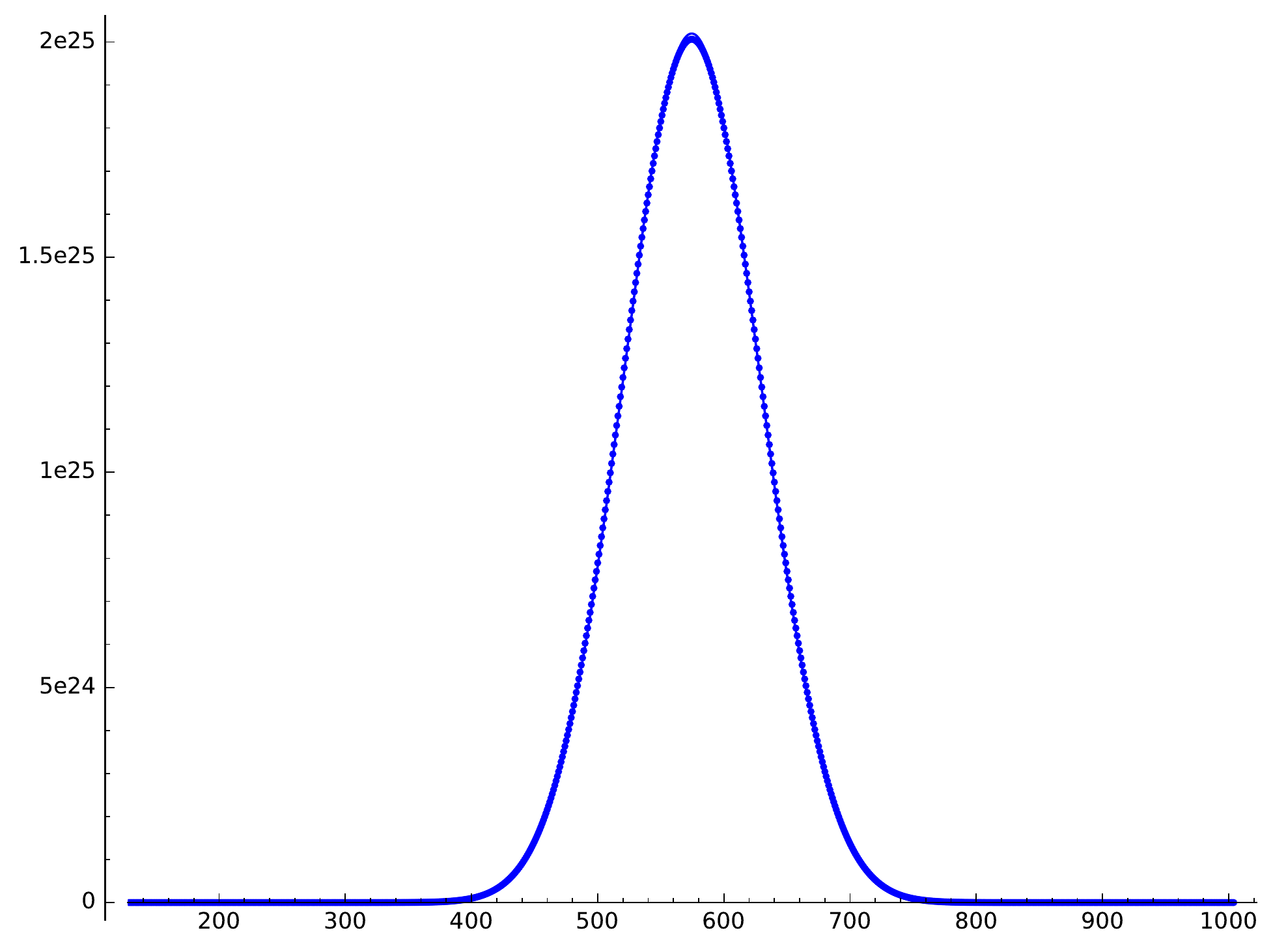}
    \caption{$\lambda = (8, 8, 7, 6, 5, 5, 5, 2, 2)$, $\aft(\lambda) = 39$}
    \label{fig:an_examples_c}
  \end{subfigure}
  \caption{Plots of $\#\{T \in \SYT(\lambda) : \maj(T) = k\}$
     as a function of $k$ for three partitions $\lambda$, overlaid
     with scaled Gaussian approximations using the same
     mean and variance.}
  \label{fig:an_examples}
\end{figure}

Let $\SYT(\lambda)$ denote the set of all standard Young tableaux of
partition shape $\lambda$.  We say $i$ is a \emph{descent} in a
standard tableau $T$ if $i+1$ comes before $i$ in the row reading word
of $T$, read from bottom to top along rows in English notation.
Equivalently, $i$ is a descent in $T$ if $i+1$ appears in a lower row
in $T$.  Let $\maj(T)$ denote the \emph{major index statistic} on
$\SYT(\lambda)$, which is again defined to be the sum of the descents
of $T$. \Cref{fig:an_examples} shows some sample distributions for the
major index on standard tableaux for three particular partition shapes.
Note that Gaussian approximations fit the data well.

In \Cref{thm:maj_Sn_an}, we simply let $n \to \infty$. For partitions,
the shape $\lambda$ may ``go to infinity'' in many different ways.
The following statistic on partitions overcomes this difficulty.

\begin{Definition}
  Suppose $\lambda$ is a partition. Let the \textit{aft} of $\lambda$ be
    \[ \aft(\lambda) \coloneqq	 |\lambda| - \max\{\lambda_1, \lambda_1'\}. \]
\end{Definition}
\noindent Intuitively, if the first row of $\lambda$ is at least as long as the
first column, then $\aft(\lambda)$ is the number of cells \textit{not}
in the first row. This definition is strongly reminiscent of a
\textit{representation stability} result of Church and Farb
\cite[Thm.~7.1]{MR3084430}, which is proved with an analysis of
the major index on standard tableaux.

Our first main result gives the analogue of \Cref{thm:maj_Sn_an} for
$\maj$ on $\SYT(\lambda)$. In particular, it completely classifies which
sequences of partition shapes give rise to asymptotically normal
sequences of $\maj$ statistics on standard tableaux.

\begin{Theorem}\label{thm:an}
  Suppose $\lambda^{(1)}, \lambda^{(2)}, \ldots$ is a sequence of
  partitions, and let $\cX_N = \cX_{\lambda^{(N)}}[\maj]$ be the corresponding
  random variables for the $\maj$ statistic on $\SYT(\lambda^{(N)})$.
  Then, the sequence $\cX_1, \cX_2, \ldots$ is asymptotically normal if and only if
  $\aft(\lambda^{(N)}) \to \infty$ as $N \to \infty$.
\end{Theorem}

\begin{Remark}
  In \Cref{sec:SYT_diag_an}, we more generally consider $\maj$
  on $\SYT(\underline{\lambda})$ where $\underline{\lambda}$ is a
  \textit{block diagonal} skew partition. See \cite[\S2]{BKS-zeroes} for
  further representation-theoretic motivation
  and \cite[Thm.~6.3]{BKS-zeroes} for the classification of the support
  of $\maj$ on $\SYT(\underline{\lambda})$.

  The generalization of \Cref{thm:an} to $\SYT(\underline{\lambda})$
  is \Cref{thm:an_diag}. Special cases of \Cref{thm:an_diag} include
  Canfield--Janson--Zeilberger's main result in \cite{MR2794017}
  classifying asymptotic normality for $\inv$ or $\maj$ on words
  (though see \cite{MR2925927c} for earlier, essentially equivalent
  results due to Diaconis \cite{MR964069}).
  The case of words generalizes
  \Cref{thm:maj_Sn_an}. The $\lambda^{(N)} = (N, N)$ case of
  \Cref{thm:an} also recovers the main result of Chen--Wang--Wang
  \cite{Chen-Wang-Wang.2008}, giving asymptotic normality for
  $q$-Catalan coefficients.
\end{Remark}

Our proof of \Cref{thm:an} relies on the \textit{method of moments},
which requires useful descriptions of the moments of $\cX_\lambda[\maj]$.
Adin--Roichman \cite{MR1841639} gave exact formulas for the mean
and variance of $\cX_\lambda[\maj]$ in terms of the
\textit{hook lengths} of $\lambda$. Their argument leverages the
following $q$-analogue of the celebrated Frame--Robinson--Thrall
Hook Length Formula \cite[Thm. 1]{Frame-Robinson-Thrall.1954}
(obtained by setting $q=1$):
\begin{equation}\label{eq:SYT_maj_hook}
  \SYT(\lambda)^{\maj}(q)
    \coloneqq \sum_{T \in \SYT(\lambda)} q^{\maj(T)}
    = q^{\rank(\lambda)} \frac{[n]_q!}{\prod_{c \in \lambda} [h_c]_q},
\end{equation}
where $h_c$ denotes the hook length of a cell $c$ in $\lambda$
and
$\rank(\lambda) \coloneqq \sum_{i \geq 1} (i-1) \lambda_i$.
Equation~\eqref{eq:SYT_maj_hook} is due to Stanley
\cite[Cor.~7.21.5]{ec2} and is
strongly related to the stable principal specialization of Schur functions
by the identity $s_\lambda(1, q, q^2, \ldots) =
\SYT(\lambda)^{\maj}(q)/\prod_{i=1}^{|\lambda|} (1-q^i)$
\cite[Prop.~7.19.11]{ec2}.

In fact, formulas for the $d$th moment $\mu_d^\lambda$,
$d$th central moment $\alpha_d^\lambda$, and
$d$th \textit{cumulant} $\kappa_d^\lambda$ of $\maj$ on
$\SYT(\lambda)$ may be
derived
from
\eqref{eq:SYT_maj_hook}. The most elegant of these
formulas is for the cumulants, from which the moments and
central moments are all easy to compute.
\begin{Theorem}\label{thm:SYT_moments}
  Let $\lambda \vdash n$ and $d \in \bZ_{> 1}$. We have
  \begin{equation}\label{eq:SYT_cumulants}
    \kappa_d^\lambda =  \frac{B_d}{d}
    \left[ \sum_{j=1}^n j^d - \sum_{c \in \lambda} h_c^d \right]
  \end{equation}
  where $B_0, B_1, B_2, \ldots = 1, \frac{1}{2}, \frac{1}{6}, 0, -\frac{1}{30}, 0, \frac{1}{42},
  0, \ldots$ are the Bernoulli numbers.
\end{Theorem}
\noindent See \Cref{thm:prod_cumulants} for a
generalization of \eqref{eq:SYT_cumulants}
along with exact formulas for the moments and
central moments. See \Cref{rem:prod_cumulants} for
the some of the history of this formula.

\begin{Remark}
  For ``most'' partition shapes, one expects the term $\sum_{j=1}^n j^d$
  in \eqref{eq:SYT_cumulants} to dominate $\sum_{c \in \lambda} h_c^d$,
  in which case asymptotic normality is quite straightforward. However,
  for some shapes there is a very large amount of cancellation in
  \eqref{eq:SYT_cumulants} and determining the limit law can be quite
  subtle.

  While $\cX_{\lambda}[\maj]$ can be written as
  the sum of scaled indicator random variables $D_1, 2D_2,3D_3, \ldots,$ $(n-1)D_{n-1}$
  where $D_i$ determines if there is a descent at position $i$,
  the $D_i$ are not at all independent,
  so one may not simply apply standard central limit theorems.
  Interestingly, the $D_i$ are identically distributed \cite[Prop.~7.19.9]{ec2}.
  The lack of independence of the $D_i$'s likewise complicates related work
  by Fulman \cite{MR1652841} and Kim--Lee \cite{1803.10457}
  considering the limiting distribution of descents in certain classes of
  permutations.
\end{Remark}

The non-normal continuous limit laws for $\maj$ on $\SYT(\lambda)$ turn out
to be the \textit{Irwin--Hall distributions} $\cIH_M \coloneqq \sum_{k=1}^M
\cU[0, 1]$, which are the sum of $M$ i.i.d.~continuous $[0, 1]$ random variables.
The following result completely classifies all possible limit laws for
$\maj$ on $\SYT(\lambda)$ for any sequence of partition shapes. See
\Cref{thm:all_limits_diag} for the generalization to block diagonal skew shapes.

\begin{Theorem}\label{thm:all_limits}
  Let $\lambda^{(1)}, \lambda^{(2)}, \ldots$ be a sequence of partitions.
  Then $(\cX_{\lambda^{(N)}}[\maj]^*)$ converges in distribution if and only if
  \begin{enumerate}[(i)]
    \item $\aft(\lambda^{(N)}) \to \infty$; or
    \item $|\lambda^{(N)}| \to \infty$ and $\aft(\lambda^{(N)}) \to M < \infty$; or
    \item the distribution of $\cX_{\lambda^{(N)}}^*[\maj]$ is eventually constant.
  \end{enumerate}
  The limit law is $\cN$ in case (i), $\cIH_M^*$ in
  case (ii), and discrete in case (iii).
\end{Theorem}

Case (iii) naturally leads to the question, when does $\cX_{\lambda}^*[\maj]
= \cX_{\mu}^*[\maj]$? Such a description in terms of hook lengths is
given in \Cref{thm:discrete_distributions}. \Cref{thm:all_limits} naturally
raises several open questions and conjectures concerning unimodality,
log-concavity, and local limit theorems, which are described in \Cref{sec:future}.

\begin{Example}
  We illustrate each possible limit in \Cref{thm:all_limits}. For (i), let
  $\lambda^{(N)} \coloneqq (N, \lfloor\ln N\rfloor)$, so that
  $\aft(\lambda^{(N)}) = \lfloor\ln N\rfloor \to \infty$ and the
  distributions are asymptotically normal.
  For (ii), fix $M \in \bZ_{\geq 0}$ and let $\lambda^{(N)}
  \coloneqq (N+M, M)$, so that $\aft(\lambda^{(N)}) = M$ is constant
  and the distributions converge to $\Sigma_M^*$.
  For (iii), let $\lambda^{(2N)} \coloneqq (12, 12, 3, 3, 3, 2, 2, 1, 1)$
  and $\lambda^{(2N+1)} \coloneqq (15, 6, 6, 6, 4, 2)$, which have the
  same multisets of hook lengths despite not being transposes of each other, and
  consequently the same normalized $\maj$ distributions.
\end{Example}

The rest of the paper is organized as follows. In
\Cref{sec:background}, we give background focused on cumulants aimed
at the combinatorial audience. In \Cref{sec:back_comb}, we collect
combinatorial background on permutations, tableaux, etc, aimed
more at the probabilistic audience.  In \Cref{sec:baj_inv}, we analyze
$\baj - \inv$ on $S_n$ as an introductory example.  In
\Cref{sec:SYT_diag_an}, we classify when $\maj$ on
$\SYT(\underline{\lambda})$ is asymptotically normal. In
\Cref{sec:SYT_diag_IH}, we determine the remaining continuous limit
laws for $\maj$ on $\SYT(\underline{\lambda})$.  In
\Cref{sec:SYT_discrete}, we characterize the possible discrete
distributions for $\maj$ on $\SYT(\lambda)$ in terms of hook
lengths. Finally, \Cref{sec:future} lists conjectures concerning
unimodality, log-concavity, and local limit theorems.

\section{Background on cumulants}\label{sec:background}

In this section, we review some standard terminology and results on
generating functions, random variables, and asymptotic normality,
with a focus on \textit{cumulants}. An excellent source for many
further details in this area can be found in Canfield's
Chapter 3 of \cite{MR3408702}.

\subsection{Exponential generating functions}\label{ssec:egf}

% Topics:
% OGF's and EGF's
% Composition of EGF's, moment to cumulant version
% Polynomial equivalence and recursions

We now introduce our notation for exponential generating functions
and the
Bernoulli numbers, which will be used with cumulants shortly.

\begin{Definition}
  Given a rational sequence $(g_d)_{d=0}^\infty = (g_0, g_1, \ldots)$,
  the corresponding \textit{ordinary generating function} is
    \[ O_g(t) \coloneqq \sum_{d \geq 0} g_d t^d \]
  and the corresponding \textit{exponential generating
  function} is
    \[ E_g(t) \coloneqq \sum_{d \geq 0} g_d \frac{t^d}{d!}. \]
  Conversely, any rational power series
    \[ F(t) = \sum_{d \geq 0}  \ f_d t^d = \sum_{d \geq 0} d!  f_d
        \frac{t^d}{d!} \]
  is the ordinary generating function of the sequence
  $(f_d)_{d=0}^\infty = (f_0, f_1, \ldots)$ and the exponential
  generating function of the sequence $(d! f_d)_{d=0}^\infty$.
  The exponential generating functions we will encounter will all
  have a positive radius of convergence.
\end{Definition}

It is easy to describe products, quotients and compositions of
generating functions.  We recall in particular a formula for compositions of
exponential generating functions for later use. Given two
rational sequences $f=(f_d)_{d=0}^\infty$, $g=(g_d)_{d=0}^\infty$ such
that $f_0=0$ and $g_0=1$, the composition of their exponential generating
functions $E_g \circ E_f$ is again an exponential generating function for
a rational sequence $h$, say $E_h(t) = E_g (E_f(t))$.  For example, if
$E_f(t)= \sum f_d t^d/d!$ and  $E_g(t)=e^t$, so $g_i=1$ for all $i$, then by
\cite[Cor.~5.1.6]{ec2}, the corresponding sequence
$(h_d)_{d=0}^\infty$ is given by $h_0 =1$ and, for $d \geq 1$,
\begin{equation}\label{eq:egf_comp}
  h_d = \sum_{\pi \in \Pi_d}  \prod_{b \in \pi} f_{|b|},
\end{equation}
where $\Pi_d$ is the collection of all set partitions
$\pi=\{b_1,b_2,\ldots, b_k\}$ of $\{1,2,\ldots, d\}$. Collecting
together $S_d$-orbits of $\Pi_d$ in \eqref{eq:egf_comp} quickly gives
\begin{equation}\label{eq:egf_comp_part}
  h_d = \sum_{\lambda \vdash d} \frac{d!}{z_\lambda} \prod_i
    \frac{f_{\lambda_i}}{(\lambda_i-1)!}
\end{equation}
where if $\lambda$ has $m_i$ parts of length $i$, then $z_\lambda
\coloneqq 1^{m_1} 2^{m_2} \cdots m_1! m_2! \cdots$. A more computationally
efficient, recursive approach to \eqref{eq:egf_comp} is the formula
\cite[Prop.~5.1.7]{ec2}
\begin{equation}\label{eq:egf_comp_rec}
  h_d  = f_d + \sum_{m=1}^{d-1} \binom{d-1}{m-1} f_m h_{d-m}.
\end{equation}

\begin{Example}
  The \emph{Bernoulli numbers} $(B_d)_{d=0}^\infty$ are rational
  numbers determined by the exponential generating function
  $E_B(t)\coloneqq t/(1-e^{-t})$.  The first few terms in the sequence are
   \[ B_0=1,\  B_1=\frac{1}{2},\  B_2=\frac{1}{6},\  B_3=0,\  B_4=-\frac{1}{30},\
       B_5=0,\  B_6=\frac{1}{42},\ \]
    \[ B_7=0,\  B_8=-\frac{1}{30},\   B_9=0,\  B_{10}=\frac{5}{66},\ B_{11}=0,\
        B_{12}=-\frac{691}{2730}. \]
  The \emph{divided Bernoulli numbers} are
  given by $\frac{B_d}{d}$ for $d \geq 1$. Their exponential generating
  function $E_D(t)$ satisfies $1 + t \frac{d}{dt} E_D(t) = E_B(t)$, from
  which it follows that
    \[ E_D(t) \coloneqq \sum_{d \geq 1} \frac{B_d}{d} \frac{t^d}{d!}
        = \log\left(\frac{e^t - 1}{t}\right). \]
      We caution that a common alternate convention
for Bernoulli numbers
      uses
  $B_1 = -\frac{1}{2}$ with all other entries the same,
  corresponding with the exponential generating
  function $t/(e^t - 1)$.
\end{Example}

The Bernoulli numbers have many interesting properties; see
\cite{mazur, wiki:bernoulli} and \cite[Section 6.5]{GKP}.  For example,
they appear in the polynomial expansion of the sums of $d$th powers,
\begin{equation}\label{eq.Bernoulli.sum.expansion}
  \sum_{k=1}^{n} k^d = \frac{1}{d+1}\sum_{k=0}^{d}
  \binom{d+1}{k} B_k \ n^{d+1-k}.
\end{equation}
  Compare the formula for sums of $d$th powers to the Riemann zeta
  function $\zeta(s)=\sum_{n=1}^{\infty} \frac{1}{n^s}$ which can be
  evaluated at complex values $s \neq 1$ by analytic
  continuation. The divided Bernoulli numbers which appear in our
  formula \eqref{eq:SYT_cumulants} satisfy
  $\frac{B_d}{d} = -\zeta(1-d)$.

\subsection{Probabilistic generating functions}\label{ssec:pgfs}

% Topics:
% PDF, PMF, CDF, alpha's, mu's
% MGF, CF, connection to PGF and OGF
% Cumulants
% 5 nice properties of cumulants
% alpha, mu, kappa polynomial equivalence, recursion

We next review basic vocabulary and notation for moments and cumulants
of random variables. All random variables we encounter will have
moments of all orders. See \cite{MR1324786} for more details.

\begin{Definition}
  Let $\cX$ be a real-valued random variable where either $\cX$ is
  continuous with probability density function
  $f \colon \bR \to \bR_{\geq 0}$ or $\cX$ is discrete with probability
  mass function $f \colon \bZ \to \bR_{\geq 0}$. The
  \textit{cumulative distribution function} (CDF) of $\cX$ is given by
    \[ F(t) \coloneqq \int_{-\infty}^t f(x)\,dx
        \qquad\text{ or }\qquad
        F(t) \coloneqq \sum_{k \leq t} f(k) \]
  depending on whether $\cX$ is continuous or discrete.
  For any continuous real-valued function $g$, there is an associated
  random variable $g(\cX)$. The \emph{expectation} of $g(\cX)$ is given by
    \[ \bE[g(\cX)] \coloneqq \int_{\bR} g(x) f(x)\,dx
        \qquad\text{ or }\qquad
        \bE[g(\cX)] \coloneqq \sum_{k=-\infty}^\infty g(k) f(k). \]
  The \textit{mean} and \textit{variance} of $\cX$ are, respectively,
    \[ \mu \coloneqq \bE[\cX]
        \qquad\text{ and }\qquad
        \sigma^2 \coloneqq \bE[(\cX-\mu)^2]. \]
  For $d \in \bZ_{\geq 0}$, the \textit{$d$th moment} and
  \textit{$d$th central moment} of $\cX$ are, respectively,
    \[ \mu_d \coloneqq \bE[\cX^d]
        \qquad\text{ and }\qquad \alpha_d \coloneqq \bE[(\cX-\mu)^d]. \]
  The \textit{moment-generating function} of $\cX$ is
    \[ M_{\cX}(t)\coloneqq\bE[e^{t\cX}] = \sum_{d=0}^\infty \mu_d \frac{t^d}{d!}, \]
  which for us will always have a positive radius of convergence.
  The \textit{characteristic function} of $\cX$ is
      \[ \phi_{\cX}(t)\coloneqq\bE[e^{it\cX}], \]
  which exists for all $t \in \bR$ and which is the Fourier transform
  of $f$, the density or mass function associated to $\cX$.
\end{Definition}

\begin{Example}\label{ex:pgf_mgf_cf}
  Let $W$ be a finite set with an integer statistic
  $\stat \colon W \to \bZ_{\geq 0}$.  We will use the notation
  \[ W^{\stat}(q) \coloneqq \sum_{w \in W}
    q^{\stat(w)} \] for the corresponding polynomial generating
  function.  If $ W^{\stat}(q) = \sum c_k q^k$, define a random
  variable $\cX$ associated with $\stat \colon W \to \bZ_{\geq 0}$
  sampled uniformly on $W$ by $\bP(\cX =k) = c_k/\#W.$ The
  \textit{probability generating function} for $\cX$ is
    \[ \bE[q^{\cX}] = \frac{1}{\#W} W^{\stat}(q)
        \coloneqq \frac{1}{\#W} \sum_{w \in W} q^{\stat(w)}. \]
   Letting $q=e^t$, an easy computation shows that the moment-generating
  function and characteristic function of $\cX$ are
    \[ M_{\cX}(t) = \frac{1}{\#W} W^{\stat}(e^t)
        \qquad\text{ and }\qquad
        \phi_{\cX}(t) = \frac{1}{\#W} W^{\stat}(e^{it}). \]
  These expressions reveal an intimate connection between the study
  of generating functions of combinatorial statistics evaluated on the
  unit circle and the underlying probability distribution via the
  Laplace and Fourier transforms.  In particular, the distribution determines the
  characteristic function and the moment-generating function, and
  conversely each of these determines the distribution.
\end{Example}

\begin{Definition}
 The
  \textit{cumulants} $\kappa_1, \kappa_2, \ldots$ of $\cX$ are defined
  to be the coefficients of the exponential generating function
    \[ K_{\cX}(t) \coloneqq \sum_{d=1}^\infty \kappa_d \frac{t^d}{d!}
       \coloneqq \log M_{\cX}(t) = \log \bE[e^{t\cX}]. \]
\end{Definition}

While cumulants of random variables may initially be less intuitive than
moments, they lead to nicer formulas in many cases, including
\Cref{thm:SYT_moments}, and they often have more useful properties.
See \cite{Novak-Sniady.2011} for some history and applications.
We will use the following properties of cumulants. The proofs are
straightforward from the definitions.

\medskip

\begin{enumerate}[1.]
  \item \textit{(Familiar Values)} The first three cumulants are
    $\kappa_1 = \mu$, $\kappa_2 = \sigma^2$, and
    $\kappa_3 = \alpha_3$.  The higher cumulants typically differ from
    the moments and central moments.
  \medskip
  \item \textit{(Shift Invariance)} The second and higher cumulants of $\cX$ agree
    with those for $\cX-c$ for  $c \in \bR$.
  \medskip
  \item \textit{(Homogeneity)} The $d$th cumulant of $c\cX$ is
    $c^d\kappa_d$ for $c \in \bR$.
  \medskip
  \item \textit{(Additivity)} The cumulants of the sum of
    \textit{independent} random variables are the sums of the cumulants.
  \medskip
  \item \textit{(Polynomial Equivalence)} The cumulants, moments, and central
    moments are determined by polynomials in any one of these three sequences.
\end{enumerate}
\medskip

The polynomial equivalence property can be made explicit by the
results in \Cref{ssec:egf}. Equation \eqref{eq:egf_comp_rec}
allows us to express the $d$th moment of $\cX$ as a polynomial function
of the first $d$ cumulants of $\cX$ and vice versa via the recurrence
\begin{align}\label{eq:moment_to_cumulant}
  \mu_d = \kappa_d + \sum_{m=1}^{d-1}
    \binom{d-1}{m-1} \kappa_m \mu_{d-m}.
\end{align}
Using the shift invariance property of cumulants, the corresponding
formula for the central moments in terms of the cumulants can be
obtained from \eqref{eq:moment_to_cumulant} by setting
$\kappa_1=0$ and leaving the other cumulants alone. This gives,
for $d>1$,
\begin{align}\label{eq:central_moment_to_cumulant}
  \alpha_d = \kappa_d + \sum_{m=2}^{d-2}
    \binom{d-1}{m-1}\kappa_m \alpha_{d-m}.
\end{align}
For instance, at $d=3$ we have
  \[ \mu_3 =   \kappa_3 + 3\kappa_2 \kappa_1 + \kappa_1^3. \]
Setting $\kappa_1 = 0$ yields $\alpha_3 = \kappa_3$ as mentioned
above.

\subsection{Cumulant formulas}\label{ssec:cumulants}

% Topics:
% N cumulants
% U[0, 1] cumulants
% Uniform discrete cumulants
% Cumulants of quotients of products of q-integers (HZ generally)
% Adin--Roichman mean, variance formulas
% alpha_d^\lambda, \mu_d^\lambda, \kappa_d^\lambda formulas (maybe diag case?)

Next we describe the cumulants of some well-known distributions and
use one of them to deduce a result of Hwang--Zacharovas, which
immediately yields \Cref{thm:SYT_moments} as a corollary.

\begin{Example}
  Let $\cX = \cN(\mu, \sigma^2)$ be the normal random variable with
  mean $\mu$ and variance $\sigma^2$. The density function of $\cX$
  is $f(x; \mu, \sigma^2) = \frac{1}{\sigma \sqrt{2\pi}}
  \exp\left(-\frac{(x - \mu)^2}{2\sigma^2}\right)$.
  Taking the Fourier transform gives the characteristic function
  $\bE[e^{it\cX}] = \exp\left(i\mu t - \frac{1}{2}\sigma^2 t^2\right)$,
  so the moment-generating function is
  $\bE[e^{t\cX}] = \exp\left(\mu t + \frac{1}{2} \sigma^2 t^2\right)$
  and the cumulants are
  \begin{align}\label{eq:normal_cumulants}
    \kappa_d
       = \begin{cases}
            \mu & d=1, \\
            \sigma^2 & d=2, \\
            0 & d \geq 3.
          \end{cases}
  \end{align}
  Using \eqref{eq:egf_comp_part} to compute the central moments of $\cX$ from
  \eqref{eq:normal_cumulants}, we effectively set $\kappa_1=0$ and note
  that only $\lambda=(2, 2, \ldots, 2) = (2^{d/2})$ contributes, in which case
  $\alpha_d = \kappa_2^{d/2} d!/(2^{d/2} (d/2)!)$. It follows that
    \[ \alpha_d =
        \begin{cases}
          0 & \text{if $d$ is odd}, \\
          \sigma^d (d-1)!! & \text{if $d$ is even}.
        \end{cases} \]
\end{Example}

\begin{Example}\label{ex:Ucont}
  Let $\cU = \cU[0, 1]$ be the continuous uniform random variable whose
  density takes the value $1$ on the interval $[0,1]$ and $0$ otherwise.
  Then the moment generating function is
  $M_{\cU}(t) = \int_{0}^1 e^{tx}dx = (e^{t} - 1)/t$, so the
  cumulant generating function $\log M_{\cU}(t)$ coincides with
  the exponential generating function for the divided Bernoulli numbers
  from \Cref{ssec:egf}. That is, $\kappa_d^{\cU} = B_d/d$ for $d \geq 1$.

  Recall from \Cref{sec:intro}, $\cIH_m$ is the \textit{Irwin--Hall}
  distribution obtained by adding $m$ independent, identically
  distributed $\cU[0, 1]$ random variables.  By Additivity, the $d$th
  cumulant of $\cIH_m$ is $mB_d/d$.  More generally, let
  $\cS \coloneqq \sum_{k=1}^m \cU[\alpha_k, \beta_k]$ be the sum of
  $m$ independent uniform continuous random variables. Then the $d$th
  cumulant of $\cS$ for $d \geq 2$ is
  \begin{equation}\label{eq:Ucont.cumulants}
    \kappa_d^{\cS} = \frac{B_d}{d} \sum_{k=1}^m (\beta_k - \alpha_k)^d
    \end{equation}
    by the Homogeneity and Additivity Properties of cumulants.
\end{Example}

\begin{Example}\label{ex:Un}
  Let $\cU_n$ be the discrete uniform random variable supported
  on $\{0, 1, \ldots, n-1\}$. The probability generating function for $\cU_n$ is
  $[n]_q/n \coloneqq (q^n-1)/(n(q-1))$, so the cumulant generating function is
    \[ \log M_{\cU_n}(t) = \log\left(\frac{e^{nt} - 1}{n(e^t - 1)}\right)
        = \log\left(\frac{e^{nt} - 1}{nt}\right)
           - \log\left(\frac{e^t - 1}{t}\right). \]
  It follows that for $d \geq 1$, the divided Bernoulli numbers arise
  again in this context,
  \begin{equation}\label{eq:cumulants_Un}
    \kappa_d^{\cU_n} = \frac{B_d}{d} (n^d - 1).
  \end{equation}
\end{Example}

Product formulas for polynomials such as Stanley's
formula \eqref{eq:SYT_maj_hook}
give rise to explicit formulas for cumulants
and moments according to the following theorem. The proof is
immediate from \Cref{ex:Un} and the exponential generating
function identity \eqref{eq:egf_comp_part}.

\begin{Theorem}\label{thm:prod_cumulants}
  Suppose $\{a_1, \ldots, a_m\}$ and $\{b_1, \ldots, b_m\}$ are
  multisets of positive integers such that
  \[ P(q) = \frac{\prod_{k=1}^m [a_k]_q}{\prod_{k=1}^m [b_k]_q}
    =\sum c_k q^k
    \in
    \bZ_{\geq 0}[q], \]
  so in particular each $c_k \in \bZ_{\geq 0}$.
  Let $\cX$ be a discrete random variable with
  $\bP[\cX=k]=c_k/P(1)$. Then the $d$th
  cumulant of $\cX$ is
  \begin{equation}\label{eq:prod_cumulants}
    \kappa_d^{\cX} =  \frac{B_d}{d} \sum_{k=1}^m (a_k^d - b_k^d)
  \end{equation}
  where $B_d$ is the $d$th Bernoulli number (with $B_1 =
  \frac{1}{2}$).
  Moreover, the $d$th central moment of $\cX$ is
  \begin{equation}\label{eq:prod_central_moments}
    \alpha_d = \sum_{\substack{\lambda \vdash d \\ \text{has all parts even}}}
      \frac{d!}{z_\lambda} \prod_{i=1}^{\ell(\lambda)}
        \frac{B_{\lambda_i}}{\lambda_i!}
        \left[\sum_{k=1}^m \left(a_k^d - b_k^d\right)\right].
  \end{equation}
  and the $d$th moment of $\cX$ is
  \begin{equation}\label{eq:prod_moments}
    \mu_d = \sum_{\substack{\lambda \vdash d \\ \text{has all parts either}
      \\ \text{even or size $1$}}}
      \frac{d!}{z_\lambda} \prod_{i=1}^{\ell(\lambda)}
        \frac{B_{\lambda_i}}{\lambda_i!}
        \left[\sum_{k=1}^m \left(a_k^d - b_k^d\right)\right].
  \end{equation}
\end{Theorem}

\begin{Remark}\label{rem:prod_cumulants}
  Equation \eqref{eq:prod_cumulants} appeared explicitly in the work of
  Hwang--Zacharovas \cite[\S4.1]{MR3346464} building on the work
  of Chen--Wang--Wang \cite[Thm. 3.1]{Chen-Wang-Wang.2008},
  who in turn used an argument going back at least to Sachkov
  \cite[\S1.3.1]{MR1453118}. It was rediscovered
  experimentally through \eqref{eq:prod_moments} by the present
  authors, and by Thiel--Williams \cite{Theil-Williams-2018}.
\end{Remark}

One frequently encounters polynomials of the form $q^\beta P(q)$
for some $\beta \in \bZ_{\geq 0}$, as in \eqref{eq:SYT_maj_hook}.
The formulas in \Cref{thm:prod_cumulants} remain valid in this
case except
that one must add $\beta$ to the expression for $\kappa_1$ and
add $\beta$ to each factor in the product in
\eqref{eq:prod_moments} for which $\lambda_i = 1$.

\begin{Remark}
  The generating function machinery used to construct the cumulants in
  \eqref{eq:prod_cumulants} works whether or not the function $P(q)$
  is polynomial.  The corresponding $\kappa_d$'s are called
  \emph{formal cumulants} in the literature.
\end{Remark}

\subsection{Asymptotic normality}\label{ssec:an}

% Topics:
% Define convergence in distribution, perhaps mention weak convergence
% Define asymptotic normality
% de Moivre-Laplace example
% CJZ on words example
% Chen--Wang--Wang q-Catalan example?
% AN table

Asymptotic normality is a very old topic lying at the intersection
of probability and combinatorics. For an introduction, we
recommend Canfield's Chapter 3 in \cite{MR3408702}.

\begin{Definition}
  Let $\cX_1, \cX_2, \ldots$ and $\cX$ be real-valued random variables
  with cumulative distribution functions $F_1, F_2, \ldots$ and $F$,
  respectively. We say $\cX_1, \cX_2, \ldots$ \textit{converges in distribution}
  to $\cX$, written $\cX_n \Rightarrow \cX$, if for all $t \in \bR$ at
  which $F$ is continuous we have
    \[ \lim_{n \to \infty} F_n(t) = F(t). \]
\end{Definition}

Recall from the introduction that for a real-valued random variable
$\cX$ with mean $\mu$ and variance $\sigma^2 > 0$, the
corresponding \textit{normalized random variable} is
  \[ \cX^* \coloneqq \frac{\cX-\mu}{\sigma}. \]
Observe that $\cX^*$ has mean $\mu^* = 0$ and variance
${\sigma^*}^2 = 1$. The moments and central moments of $\cX^*$
agree for $d \geq 2$ and are given by
  \[\mu_d^* = \alpha_d^* = \alpha_d/\sigma^d. \]
Similarly, the cumulants of $\cX^*$ are given by $\kappa_1^* = 0$,
$\kappa_2^* = 1$, and $\kappa_d^* = \kappa_d/\sigma^d$ for
$d \geq 2$.

\begin{Definition}
  Let $\cX_1, \cX_2, \ldots$ be a sequence of real-valued random variables.
  We say the sequence is \textit{asymptotically normal} if
  $\cX_n^* \Rightarrow \cN(0, 1)$.
\end{Definition}

The ``original'' asymptotic normality result is as follows.
Let $2^{[n]}$ be the set of all subsets of
$[n] \coloneqq \{1, 2, \ldots, n\}$. Let $\cX_{2^{[n]}}[\size]$
denote the random variable given by the cardinality, where
$2^{[n]}$ is given the uniform distribution. This has the same
distribution as the number of heads after $n$ fair coin flips,
so the probability generating function up to normalization is
$(1+q)^n$. The following
result is credited to de Moivre and Laplace; see
\cite[Theorem~3.2.1]{MR3408702} for further discussion.

\begin{Theorem}[de Moivre--Laplace]\label{thm:bin_an}
  The sequence $\cX_{2^{[n]}}[\size]$ is asymptotically normal.
\end{Theorem}

Asymptotic normality results for combinatorial statistics are plentiful.
See \Cref{tab:an} for more examples and further references.

\subsection{The method of moments}\label{ssec:moments}

% Topics:
% Levy continuity (general limit laws)
% Frechet-Shohat (general limit laws)
% Cumulant version of Frechet-Shohat
% Asymptotic normality test by normalized cumulant decay condition

We next describe two standard criteria for establishing
asymptotic normality or more generally convergence in
distribution of a sequence of random variables.

\begin{Theorem}[L\'evy's Continuity Theorem, {\cite[Theorem~26.3]{MR1324786}}]\label{thm:levy}
  A sequence $\cX_1, \cX_2, \ldots$ of real-valued random variables
  converges in distribution to a real-valued random variable $\cX$
  if and only if, for all $t \in \bR$,
    \[ \lim_{n \to \infty} \bE[e^{it\cX_n}] = \bE[e^{it\cX}]. \]
\end{Theorem}

\begin{Theorem}[Frech\'et--Shohat Theorem,
  {\cite[Theorem~30.2]{MR1324786}}]\label{thm:moments}
  Let $\cX_1, \cX_2, \ldots$ be a sequence of real-valued random
  variables,
  and let $\cX$ be a real-valued random variable. Suppose
  the moments of $\cX_n$ and $\cX$ all exist and the moment
  generating functions all have a positive radius of convergence. If
  \begin{equation}\label{eq:moments_criterion}
    \lim_{n \to \infty} \mu_d^{\cX_n} = \mu_d^{\cX} \hspace{.5cm} \forall
     d \in \bZ_{\geq 1},
  \end{equation}
  then $\cX_1, \cX_2, \ldots$ converges in distribution to $\cX$.
\end{Theorem}

By \Cref{thm:levy}, we may test for asymptotic normality by
checking if the normalized characteristic functions tend point-wise to
the characteristic function of the standard normal. Likewise
by \Cref{thm:moments} we may instead perform the check on the
level of individual normalized moments, which is often referred to as the
\textit{method of moments}. By \eqref{eq:moment_to_cumulant}
we may further replace the moment condition
\eqref{eq:moments_criterion} with the cumulant condition
\begin{equation}\label{eq:cumulants_criterion}
  \lim_{n \to \infty} \kappa_d^{\cX_n} = \kappa_d^{\cX}.
\end{equation}
For instance, we have the following explicit criterion.

\begin{Corollary}\label{cor:cumulants}
  A sequence $\cX_1, \cX_2, \ldots$ of real-valued
  random variables on finite sets is asymptotically normal
  if for all $d \geq 3$ we have
  \begin{equation}\label{eq:cumulants_criterion2}
    \lim_{n \to \infty} \frac{\kappa^{\cX_n}_d}{(\sigma^{\cX_n})^d} = 0
  \end{equation}
\end{Corollary}

In fact, one may show a converse of the Frech\'et--Shohat
theorem holds for quotients as in \Cref{thm:prod_cumulants},
though we will not have need of it here.

\subsection{Local limit theorems}\label{ssec:llts}

% Topics:
% Define local limit theorems
% CJZ local limit theorem

Asymptotic normality concerns cumulative distribution functions,
so it gives estimates for the number of combinatorial objects with a
large range of statistics. However, our original motivation was to
count combinatorial objects with a given statistic. Estimates of this
latter form are frequently referred to as \textit{local limit theorems}.
Here we review two motivating examples.

The present work was partly inspired by the following local limit
theorem due to the third author with a uniform rather than normal
limit law. For $\lambda \vdash n$, let $\maj_n \colon
\SYT(\lambda) \to [n]$ be $\maj$ modulo $n$.

\begin{Theorem}\cite[Theorem~1.9]{s17}
  For $\lambda \vdash n$, let $X_\lambda[\maj_n]$ denote the
  random variable $\maj_n$ on $\SYT(\lambda)$. Suppose
  $\#\SYT(\lambda) \geq n^5$. Then, for all $k \in [n]$,
    \[ \left|\bP[X_\lambda[\maj_n] = k] - \frac{1}{n}\right|
        < \frac{1}{n^2}. \]
\end{Theorem}

Further motivation was provided by the following analogue of
\Cref{thm:cjz_an}.

\begin{Theorem}{\cite[Theorem~4.5]{MR2794017}}
  There exists a positive constant $c$ such that for every $C$,
  the following is true. Uniformly for all compositions
  $\alpha = (\alpha_1, \ldots, \alpha_m)$ such that
  $\max_i \alpha_i \leq Ce^{c s(\alpha)}$ and all integers $k$,
    \[ \bP[X_\alpha = k] = \frac{1}{\sigma\sqrt{2\pi}}
          \left(e^{-(k-\mu)^2/(2\sigma^2)}
          + O\left(\frac{1}{s(\alpha)}\right)\right), \]
  where $X_\alpha$ denotes inversions on words of type $\alpha$.
\end{Theorem}

\section{Combinatorial background}\label{sec:back_comb}

\subsection{Combinatorial background for $\baj - \inv$ on $S_n$}\label{ssec:back_comb:baj_inv}

Here we introduce the two most well-known permutation statistics,
$\inv$ and $\maj$, as well as one unusual permutation statistic,
$\baj$.

\begin{Definition}\label{def:invmajdes}
  Let $\sigma \in S_n$ be a permutation of $\{1, \ldots, n\}$.
  Set
  \begin{align*}
    \Inv(\sigma)
      &\coloneqq \{(i, j) : i < j \text{ and } \sigma(i) > \sigma(j)\}
      &(\textit{inversion set}) \\
    \inv(\sigma)
      &\coloneqq |\Inv(\sigma)|
      &(\textit{inversion number, i.e.~length}) \\
    \Des(\sigma)
      &\coloneqq \{1 \leq i \leq n-1 : \sigma(i) > \sigma(i+1)\}
      &(\textit{descent set}) \\
    \maj(\sigma)
      &\coloneqq \sum_{i \in \Des(\sigma)} i
      &(\textit{major index}).
  \end{align*}
  Following Zabrocki \cite{math/0310301}
    for the nomenclature, we also set
    \[ \baj(\sigma) \coloneqq \sum_{i \in \Des(\sigma)} i(n-i). \]
\end{Definition}

The equidistribution of $\inv$ and $\maj$ on $S_n$ is
due to MacMahon, who also first introduced $\maj$. His proof gave
the following generating function expression for
both statistics.

\begin{Theorem}[{\cite[Art.~6]{MacMahon.1913}}]\label{thm:inv_maj_cgf}
  We have
  \[ S_n^{\inv}(q) = [n]_q! \coloneqq
    \prod_{k=1}^{n-1} (1+q+q^2 +\cdots+q^k)
    = S_n^{\maj}(q). \]
      %% modified because [n]! technically not defined yet.
  \end{Theorem}
  The statistic $\baj - \inv$ appeared in the context of extended
  affine Weyl groups and Hecke algebras in the work of Iwahori and
  Matsumoto in 1965 \cite{Iw.M}.  It is the Coxeter length function
  restricted to coset representatives of the extended affine Weyl
  group of type $A_{n-1}$ mod translations by coroots.  Stembridge and
  Waugh \cite[Remarks 1.5 and 2.3]{MR1692145} give a careful overview
  of this topic and further results.  In particular, they prove the
  following factorization formula for the generating function
  associated to $\baj - \inv$ on $S_n$.  From this factorization, the
  corresponding cumulants can be read off from
  \Cref{thm:prod_cumulants}.

\begin{Theorem}\cite{Iw.M,MR1692145}
  We have
  \begin{equation}\label{eq:baj_inv_prod}
    S_n^{\baj - \inv}(q)
      \coloneqq \sum_{\sigma \in S_n} q^{\baj(\sigma) - \inv(\sigma)}
      = n \prod_{i=1}^{n-1} \frac{[i(n-i)]_q}{[i]_q}.
  \end{equation}
\end{Theorem}

\begin{Corollary}\label{cor:baj_inv_kappa}
  The $d$th cumulant $\kappa_d^n$ for $\baj - \inv$ on $S_n$
  is
    \[ \kappa_d^n = \frac{B_d}{d} \left(\sum_{i=1}^{n-1} [i(n-i)]^d - i^d\right). \]
\end{Corollary}

\begin{Remark}
  Indeed, \eqref{eq:baj_inv_prod} holds with $S_n$ replaced by
  $\{\sigma \in S_n : \sigma(n) = k\}$ for any fixed
  $k = 1, \ldots, n$ if the factor of $n$ is deleted from the
  right-hand side. See \cite{math/0310301} for a bijective proof of
  this generalization.  In addition, \cite[Thm.~1.1]{MR1692145} gives
  another generalization of the product formula
  \eqref{eq:baj_inv_prod} to all crystallographic Coxeter groups.
\end{Remark}

\subsection{Combinatorial background for $\maj$ on $\W_\alpha$
  and $\SYT(\protect\underline{\lambda})$}\label{ssec:words_tableaux}

Here we review standard com\-bi\-na\-to\-ri\-al notions related to
words, tableaux, and their major index generating functions.

\begin{Definition}
  Given a word $w = w_1w_2\cdots w_n$ with letters
  $w_i \in \bZ_{\geq 1}$, the \textit{type} of $w$ is the sequence
  $\alpha = (\alpha_1, \alpha_2, \ldots)$ where $\alpha_i$ is the
  number of times $i$ appears in $w$.  Such a sequence $\alpha$ is
  a (weak) \textit{composition} of $n$, written as
  $\alpha \vDash n$. Trailing $0$'s are often omitted when writing
  weak compositions, so $\alpha = (\alpha_1, \alpha_2, \ldots, \alpha_m)$
  for some $m$.  Note that a word of type $(1,1,\ldots,1)\vDash n$ is a
  permutation in the symmetric group $S_n$ written in one-line
  notation. Just as for permutations, the \textit{inversion number} of
  $w$ is
    \[ \inv(w) \coloneqq \#\{(i, j) : i < j, w_i > w_j\}. \]
  The \textit{descent set} of $w$ is
    \[ \Des(w) \coloneqq \{0<i <n : w_i > w_{i+1}\}, \]
  and the \textit{major index} of $w$ is
    \[ \maj(w) \coloneqq \sum_{i \in \Des(w)} i. \]
\end{Definition}

\begin{Definition}
  Let $\alpha = (\alpha_1, \ldots, \alpha_m) \vDash n$. We use the
  following standard $q$-analogues:
\[
  \begin{array}{lllll}  \vspace{.2in}
    [n]_q &\coloneqq & 1 + q + \cdots + q^{n-1} = \frac{q^n - 1}{q - 1}, &
                                                                  \hspace{.3in}
    & \text{($q$-integer)}\\ \vspace{.2in}
    {}
   [n]_q! &\coloneqq& [n]_q [n-1]_q \cdots [1]_q, &
                                                                  \hspace{.3in}
    & \text{($q$-factorial)}\\ \vspace{.2in}

    \binom{n}{k}_q
      &\coloneqq & \frac{[n]_q!}{[k]_q! [n-k]_q!} \in \bZ_{\geq 0}[q], &
                                                                  \hspace{.3in}
    & \text{($q$-binomial)}\\ \vspace{.2in}
    \binom{n}{\alpha}_q
      &\coloneqq & \frac{[n]_q!}{[\alpha_1]_q! \cdots [\alpha_m]_q!} \in
           \bZ_{\geq 0}[q]
           &
                                                                  \hspace{.3in}
    & \text{($q$-multinomial).}\\
  \end{array}
  \]
\end{Definition}

\begin{Example}
  The identity statistic on the set $W = \{0, \ldots, n-1\}$ has
  generating function $[n]_q$. The ``sum'' statistic on
  $W = \prod_{k=1}^n \{0, \ldots, k-1\}$ has generating function
  $[n]_q!$.
\end{Example}

For $\alpha \vDash n$, let $\W_\alpha$ denote the words of type
$\alpha$. MacMahon's classic result generalizing \Cref{thm:inv_maj_cgf}
in fact shows that $\maj$ and $\inv$ have the same distribution
on $\W_\alpha$.

\begin{Theorem}[{\cite[Art.~6]{MacMahon.1913}}]\label{thm:macmahon}
  For each $\alpha \vDash n$,
  \begin{align}\label{eqn:macmahon}
     \W_\alpha^{\maj}(q) = \binom{n}{\alpha}_q = \W_\alpha^{\inv}(q).
  \end{align}
\end{Theorem}

\begin{figure}[ht]
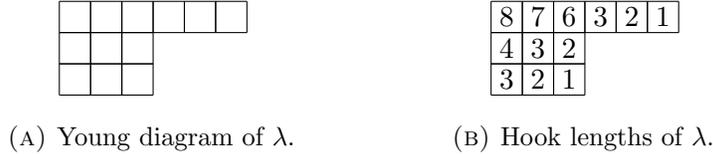

  \centering
  \begin{subfigure}[t]{0.32\textwidth}
    \centering
    \[
    \tableau{{} & {} & {} & {} & {} & {} \\ {} & {} & {} \\ {} & {} & {} }
    \]
    \caption{Young diagram of $\lambda$.}
    \label{fig:partition_a}
  \end{subfigure}
  \hspace{.2cm}
  \begin{subfigure}[t]{0.32\textwidth}
    \centering
    \[
    \tableau{{8} & {7} & {6} & {3} & {2} & {1} \\ {4} & {3} & {2} \\ {3} & {2} & {1} }
    \]
    \caption{Hook lengths of $\lambda$.}
    \label{fig:partition_c}
  \end{subfigure}
  \caption{Constructions related to the partition $\lambda=(6, 3, 3)\vdash 12$.}
  \label{fig:partition}
\end{figure}

\begin{Definition}
  A composition $\lambda \vDash n$ such that
  $\lambda_1 \geq \lambda_2 \geq \ldots$ is called a
  \textit{partition} of $n$, written as $\lambda \vdash n$. The
  \textit{size} of $\lambda$ is $|\lambda| \coloneqq n$ and the
  \textit{length} $\ell(\lambda)$ of $\lambda$ is the number of
  non-zero entries. The \textit{Young diagram} of $\lambda$ is the
  upper-left justified arrangement of unit squares called
  \textit{cells} where the $i$th row from the top has $\lambda_i$
  cells following the English notation; see \Cref{fig:partition_a}.
  The \textit{hook length} of a cell
  $c \in \lambda$ is the number $h_c$ of cells in $\lambda$ in the
  same row as $c$ to the right of $c$ and in the same column as $c$
  and below $c$, including $c$ itself; see \Cref{fig:partition_c}.  A
  \emph{corner} of $\lambda$ is any cell with hook length $1$. A
  \emph{bijective filling} of $\lambda$ is any labeling of the cells
  of $\lambda$ by the numbers $[n]=\{1,2,\ldots,n\}$.
\end{Definition}

\begin{Definition}
  A \textit{skew partition} $\lambda/\nu$ is a pair of partitions
  $(\nu, \lambda)$ such that the Young diagram of $\nu$ is contained
  in the Young diagram of $\lambda$. The cells of $\lambda/\nu$ are
  the cells in the diagram of $\lambda$ which are not in the diagram
  of $\nu$, written $c \in \lambda/\nu$.  We identify straight
  partitions $\lambda$ with skew partitions $\lambda/\varnothing$
  where $\varnothing = (0, 0, \ldots)$ is the empty partition. The
  \textit{size} of $\lambda/\nu$ is $|\lambda/\nu| \coloneqq |\lambda| - |\nu|$. The
  notions of bijective filling, hook lengths, and corners
  naturally extend to skew partitions as well.
\end{Definition}

\begin{figure}[ht]
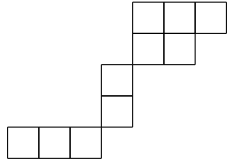

  \centering
    \[
    \tableau{& & & & {} & {} & {} \\ & & & & {} & {} \\ & & & {} \\ & & & {} \\ {} & {} & {} }
    \]
    \caption{Diagram for the skew partition
      $\lambda/\nu = 76443/4433$, which is also the block diagonal
      skew shape
      $\underline{\lambda} = ((3, 2), (1, 1),(3) )$.}
    \label{fig:diag_a}
\end{figure}

\begin{Definition}\label{def:block.diag.skew.partions}
  Given a sequence of partitions $\underline{\lambda} = (\lambda^{(1)},
  \ldots, \lambda^{(m)})$, we identify the sequence with the
  \emph{block diagonal skew partition} obtained by translating the Young diagrams of the
  $\lambda^{(i)}$ so that the rows and columns occupied by these
  components are disjoint, form a valid skew shape, and appear in order from
  top to bottom as depicted in  \Cref{fig:diag_a}.
\end{Definition}

\begin{Definition}\label{def:SYT_maj}
  A \textit{standard Young tableau} of shape $\lambda/\nu$ is a
  bijective filling of the cells of $\lambda/\nu$ such that labels
  increase to the right in rows and down columns; see
  \Cref{fig:SYTpartition}.   The set of standard Young
  tableaux of shape $\lambda/\nu$ is denoted $\SYT(\lambda/\nu)$.
  The \textit{descent set} of $T \in \SYT(\lambda/\nu)$ is the set
  $\Des(T)$ of all labels $i$ in $T$ such that $i+1$ is in a
  strictly lower row than $i$.  The \textit{major index} of $T$ is
    \[ \maj(T) \coloneqq \sum_{i \in \Des(T)} i. \]
\end{Definition}

\begin{figure}[ht]
  \centering
  \begin{subfigure}[t]{0.32\textwidth}
    \centering
    \[
    \tableau{{1} & {2} & {4} & {7} & {9} & {12} \\ {3} & {6} & {10} \\ {5} & {8} & {11} }
    \]
  \end{subfigure}
    \begin{subfigure}[t]{0.32\textwidth}
    \centering
    \[
    \tableau{& & & & & 2 & 6 \\ & & & 4 & 5 \\ 1 & 3 & 7}
    \]
  \end{subfigure}
  \caption{On the left is a standard Young tableau of straight shape
    $\lambda=(6,3,3)$ with descent set  $\{2, 4, 7, 9, 10\}$ and major
    index $32$.  On the right is a standard Young tableau of block diagonal
    skew shape $(7,5,3)/(5,3)$ corresponding to the sequence of
    partitions $\underline{\lambda} = ((2), (2), (3))$ with descent set
    $\{2,6\}$ and major index $8$. }
  \label{fig:SYTpartition}
\end{figure}

\begin{Remark}
  The block diagonal skew partitions $\underline{\lambda}$ allow us to
  simultaneously consider words and tableaux as follows.  Recall that
  $\W_\alpha$ is set of all words with type $\alpha
  = (\alpha_1, \ldots, \alpha_k)$.  Letting
  $\underline{\lambda} = ((\alpha_k), \ldots, (\alpha_1))$, we have a bijection
  \begin{equation}\label{eq:bij}
    \phi \colon \SYT(\underline{\lambda}) \too{\sim} \W_\alpha
  \end{equation}
  which sends a tableau $T$ to the word whose $i$th letter
  is the row number in which $i$ appears in $T$, counting from the \textit{bottom up}
  rather than top down.
  For example, using the skew tableau $T$ on the right of \Cref{fig:SYTpartition}, we have
  $\phi(T) = 1312231 \in \W_{(3, 2, 2)}$. It is easy to see that
  $\Des(\phi(T)) = \Des(T)$, so that $\maj(\phi(T))
  = \maj(T)$. Hence $\SYT((\alpha_1), \ldots, (\alpha_k))^{\maj}(q)
  = \W_\alpha^{\maj}(q) = \binom{n}{\alpha}_q$.
\end{Remark}

\begin{Remark}
  We also recover $q$-integers, $q$-binomials, $q$-multinomials, and
  $q$-Catalan numbers up to $q$-shifts as special cases of the
  major index generating function for tableaux given in
  \eqref{eq:SYT_maj_hook}:
  \begin{align*}
    \SYT(\lambda)^{\maj}(q)
      &= \begin{cases}
             q[n]_q & \text{if }\lambda = (n, 1), \\
             q^{\binom{k+1}{2}} \binom{n}{k}_q & \text{if }\lambda = (n-k+1, 1^k), \\
             q^n \frac{1}{[n+1]_q} \binom{2n}{n}_q & \text{if }\lambda = (n, n).
           \end{cases}
  \end{align*}
\end{Remark}

Many combinatorial statistics arise from sets indexed by more
complicated objects than the positive integers, in which case one can
``let $n \to \infty$'' in many different ways. The following result
due to Canfield, Janson, and Zeilberger illustrates a more
interesting limit.
Their result is characterized by the statistic   $s(\alpha) \coloneqq
n - m$ where $\alpha=(\alpha_1, \ldots,\alpha_\ell)\vDash n$ with $\max\{\alpha_i\}=m$.

\begin{Theorem}{\cite[Theorem~1.2]{MR2794017}}\label{thm:cjz_an}
  Let $\alpha^{(1)}, \alpha^{(2)}, \ldots$ be a sequence of compositions,
  possibly of differing lengths.
 Let $\cX_n$ be the inversion (or major index) statistic
  on words of type $\alpha^{(n)}$. Then $\cX_1, \cX_2, \ldots$ is
  asymptotically normal if and only if
    \[ s(\alpha^{(n)}) \to \infty. \]
\end{Theorem}

\begin{Remark}\label{rem:cjz_priority}
  Explorations equivalent to \Cref{thm:cjz_an} appeared significantly
  earlier than \cite{MR2794017}
  in other contexts, for instance \cite[p.~127-128]{MR964069}
  and (in the two-letter case) \cite{MR0022058}. See
  \cite{MR2925927c} for further discussion and references.
\end{Remark}

The cumulant formula for $\cX_\lambda[\maj]$,
\Cref{thm:SYT_moments}, follows immediately from
\Cref{thm:prod_cumulants} and Stanley's formula
\eqref{eq:SYT_maj_hook}. Adin and Roichman
\cite{MR1841639} had previously used
\eqref{eq:SYT_maj_hook} to compute the mean and variance of
$\cX_\lambda[\maj]$ as
  \[ \mu
      = \frac{\binom{|\lambda|}{2} - \rank(\lambda') + \rank(\lambda)}{2}
      = \rank(\lambda) + \frac{1}{2} \left[\sum_{k=1}^{|\lambda|} k
         - \sum_{c \in \lambda} h_c\right], \]
and
  \[ \sigma^2
      = \frac{1}{12} \left[\sum_{k=1}^{|\lambda|} k^2
      - \sum_{c \in \lambda} h_c^2\right]. \]

The following common generalization of Stanley's formula
\eqref{eq:SYT_maj_hook} and MacMahon's formula,
\Cref{thm:macmahon}, is well known
(e.g.\ see \cite[(5.6)]{stembridge89}).
See \cite[Thm.~2.15]{BKS-zeroes} for
other applications.

\begin{Theorem}\label{thm:diag_maj} %% block diagonal partition case
  Let $\underline{\lambda} = (\lambda^{(1)}, \ldots, \lambda^{(m)})$
  where $\lambda^{(i)} \vdash \alpha_i$ and $n = \alpha_1 + \cdots + \alpha_m$.
  Then
  \begin{equation}\label{eq:diag_maj}
    \SYT(\underline{\lambda})^{\maj}(q)
      = \binom{n}{\alpha_1, \ldots, \alpha_m}_q \cdot \prod_{i=1}^m \SYT(\lambda^{(i)})^{\maj}(q).
  \end{equation}
\end{Theorem}

\begin{Corollary}
  Let $\kappa_d^{\underline{\lambda}}$ be the $d$th cumulant of
  $\maj$ on $\SYT(\underline{\lambda})$ for $d > 1$. Then
  \begin{equation}\label{eq:cumulant_diag}
    \kappa_d^{\underline{\lambda}} = \frac{B_d}{d}
        \left(\sum_{k=1}^{|\underline{\lambda}|} k^d
        - \sum_{c \in \underline{\lambda}} h_c^d\right).
  \end{equation}
\end{Corollary}

For general skew shapes, $\SYT(\lambda/\nu)^{\maj}(q)$ does not
factor as a product of cyclotomic polynomials times $q$ to a power.
A ``$q$-Naruse'' formula due to Morales--Pak--Panova,
\cite[(3.4)]{1512.08348}, gives an analogue of
\eqref{eq:SYT_maj_hook} involving a sum over ``excited
diagrams,'' though the resulting sum has a single term
precisely for the block diagonal skew partitions $\underline{\lambda}$.

\section{Asymptotic normality for $\baj - \inv$ on $S_n$}\label{sec:baj_inv}
We give with a straightforward example which serves as
a warmup and establishes some notation. See
\Cref{ssec:back_comb:baj_inv} for background. Asymptotic normality of
$\baj - \inv$ on $S_n$ follows from the cumulant formula in \Cref{cor:baj_inv_kappa}
by the following routine calculations. Recall that $a_n \sim b_n$ means that $\lim_{n \to \infty} a_n/b_n = 1$.

\begin{Lemma}\label{lem:baj_inv.1}
  Fix $d \geq 1$. Then, as $n \to \infty$,
    \[ \sum_{i=1}^{n-1} [i(n-i)]^d - i^d
      \sim n^{2d+1} \cdot \int_0^1 x^d (1-x)^d\,dx. \]
%That is, the ratio of the left-hand side to the right-hand side tends to  $1$ as $n \to \infty$.

  \begin{proof}
    We have
    \begin{align*}
      \lim_{n \to \infty}
        \frac{\sum_{i=1}^{n-1} [i(n-i)]^d - i^d}{n^{2d+1}}
        &= \lim_{n \to \infty} \frac{1}{n} \sum_{i=1}^{n-1} \left[\left(\frac{i}{n}\right)^d
             \left(1-\frac{i}{n}\right)^d - \left(\frac{i}{n^2}\right)^d\right] \\
        &= \lim_{n \to \infty} \frac{1}{n} \sum_{i=1}^{n-1} \left(\frac{i}{n}\right)^d
             \left(1-\frac{i}{n}\right)^d \\
        &= \int_0^1 x^d (1-x)^d\,dx.
    \end{align*}
  \end{proof}
\end{Lemma}

\begin{Remark}
  The value of the integral in \Cref{lem:baj_inv.1} is well known:
  \begin{equation}\label{eq:beta}
    \int_0^1 x^d (1-x)^d\,dx = \frac{(d!)^2}{(2d+1)!}
        = \frac{1}{2d+1} \binom{2d}{d}^{-1}.
  \end{equation}
  See \cite[A002457]{oeis} for a surprisingly large number of interpretations of the
  reciprocals of these values. Equation \eqref{eq:beta} is also a very special case of the
  Selberg integral formula \cite{MR0018287}, which has many interesting connections to
  algebraic combinatorics such as those in \cite{MR3551643}.
\end{Remark}

\begin{Corollary}\label{lem:baj_inv.2}
  Fix $d \in \{1, 2, 4, 6, \ldots\}$. Let $\kappa_d^n$ be the
  $d$th cumulant of $\baj - \inv$ on $S_n$, and let
  ${\kappa_d^n}^*$ be the $d$th cumulant of the corresponding
  normalized random variable with mean $0$ and variance $1$.
  Then, uniformly for all $n$, we have
  \begin{equation} |{\kappa_d^n}^*| = \Theta(n^{1-d/2}).
  \end{equation}
  That is, there are constants $c, C>0$ depending only on $d$ such that
    \[ cn^{1-d/2} \leq |{\kappa_d^n}^*| \leq Cn^{1-d/2}. \]

  \begin{proof}
    It follows immediately from \Cref{cor:baj_inv_kappa} and
    \Cref{lem:baj_inv.1} that $|\kappa_d^n| = \Theta(n^{2d+1})$. Hence
    \[ |\kappa_d^{n*}| = |\kappa_d^n/(\kappa_2^n)^{d/2}| = \Theta(n^{2d+1 - 5d/2})
          = \Theta(n^{1-d/2}). \]
  \end{proof}
\end{Corollary}

\begin{Theorem}
  Let $\cX_n = \cX_{S_n}[\baj-\inv]$ be the random variable for the
  $\baj - \inv$ statistic taken uniformly at random from $S_n$.
  Then, $\cX_1, \cX_2, \ldots$ is asymptotically normal.

  \begin{proof}
    For fixed $d > 2$ even, we have $1-d/2 < 0$, so by
    \Cref{lem:baj_inv.2}, ${\kappa_d^n}^* \to 0$ as $n \to \infty$.
    The odd cumulants for $d > 2$ vanish
    since the odd Bernoulli numbers are $0$.
    The result now follows from \Cref{cor:cumulants}.
  \end{proof}
\end{Theorem}

\begin{Remark}
  A key step in the above argument was to show that the variance
  $\sigma_n^2$ of $\baj - \inv$ on $S_n$ satisfies
  $\sigma_n^2 = \Theta(n^5)$. Indeed, the argument
  gives $\sigma_n^2 \sim n^5/360$. The weaker observation
  that $\sum_{i=1}^{n-1} [i(n-i)]^2$ is the
  dominant contribution to $\sigma_n^2$ is essentially enough to
  deduce asymptotic normality in this case. Our analysis of
  $\maj$ on standard tableaux includes
  non-normal limits, so more precise estimates like the above will
  become absolutely necessary. A straightforward modification
  of the above argument together with \Cref{thm:inv_maj_cgf}
  also proves \Cref{thm:maj_Sn_an}.
\end{Remark}

\section{Asymptotic normality for $\maj$ on $\SYT(\protect\underline{\lambda})$}\label{sec:SYT_diag_an}

The main result of this section, \Cref{thm:an_diag}, classifies the sequences
of block diagonal skew partitions for which $\maj$ is asymptotically normal.
% See \Cref{ssec:words_tableaux} for background.
%SB: unnecessary at this point.
We begin with a series of estimates for the differences
$\sum_{k=1}^{|\lambda/\nu|} k^d - \sum_{c \in \lambda/\nu} h_c^d$,
culminating in \Cref{cor:bound.cumulants}.

\begin{Definition}
  A \textit{reverse standard Young tableau} of shape $\lambda/\nu$
  is a bijective filling of $\lambda/\nu$ which strictly
  decreases along rows and columns. The set of reverse
  standard Young tableaux of shape $\lambda/\nu$ is denoted
  $\RSYT(\lambda/\nu)$.
\end{Definition}

\begin{Lemma}\label{lem:rsyt}
  Let $\lambda/\nu \vdash n$ and $T \in \RSYT(\lambda/\nu)$. Then
  for all $c \in \lambda/\nu$,
  \begin{align}\label{eq:rsyt1}
    T_c \geq h_c.
  \end{align}
  Furthermore, for any positive integer $d$,
  \begin{align}\label{eq:rsyt2}
    \sum_{j=1}^n j^d - \sum_{c \in \lambda/\nu} h_c^d
      &= \sum_{c \in \lambda/\nu} (T_c^d - h_c^d)
        = \sum_{c \in \lambda/\nu} (T_c - h_c) {\mathbf{h}}_{d-1}(T_c, h_c),
  \end{align}
  where ${\mathbf h}_{d-1}$ denotes the complete homogeneous symmetric function.

  \begin{proof}
    For \eqref{eq:rsyt1}, equality holds at the outer corner $c$ where $T_c=1$.
    Removing $c$ and subtracting $1$ from each remaining
    entry in $T$ allows us to induct. Equation \eqref{eq:rsyt2} follows
    immediately by rearranging the terms and factoring
    $(T_c^d - h_c^d) = (T_c- h_c)\sum_{k=0}^{d-1} T_c^{d-1-k}h_c^{k}$.
  \end{proof}
\end{Lemma}

\begin{Lemma}\label{lem:small_hook_bounds}
  Let $\lambda/\nu \vdash n$ such that
  $\max_{c \in \lambda/\nu} h_c < 0.8n$. Let $d$ be
  any positive integer. Then
  \[ \frac{n^{d+1}}{26(d+1)} - 2(0.8)^d n^d
      < \sum_{j=1}^n j^d - \sum_{c \in \lambda/\nu} h_c^d
      < \frac{n^{d+1}}{d+1} + n^d. \]

  \begin{proof}
    Using Riemmann sums for $\int_0^n x^d dx$, we obtain the bounds
    \begin{equation}\label{eq:sum.bounds}
      \frac{n^{d+1}}{d+1} < \sum_{j=1}^n j^d  < \frac{n^{d+1}}{d+1} + n^d
    \end{equation}
    for all positive integers $d,n$. The upper bound in the lemma now
    follows immediately.

%     %%
%          \hl{Should we comment on proving the
%       inequalities above in light of the fact that the absolute value of the
%       Bernoulli numbers $B_k$ grow like $2 \sqrt{2\pi k} (\frac{k}{2\pi e} )^k$ ?}
% See
      %% \url{http://oeis.org/wiki/Bernoulli_numbers#Asymptotic_approximation}}.
      %% Cite that ref too in background section

    For the lower bound, label the cells of $\lambda/\nu$ by some
    $T \in \RSYT(\lambda/\nu)$. By \eqref{eq:rsyt1}, $h_c \leq T_c$,
    and by assumption we have $h_c < 0.8n$ for all $c \in \lambda/\nu$.
    Considering the tighter of these two bounds on each summand
    and using \eqref{eq:sum.bounds} again, we have
    \begin{align*}
      \sum_{c \in \lambda/\nu} h_c^d
        &< \sum_{\substack{j \in [n] \\ j < 0.8n}} j^d
           + \sum_{\substack{j \in [n] \\ j \geq 0.8n}} (0.8n)^d \\
        &< \frac{\lfloor 0.8n\rfloor^{d+1}}{d+1} + \lfloor 0.8n\rfloor^d
             + (n-\lceil 0.8n\rceil + 1)(0.8n)^d \\
        &\leq \frac{(0.8n)^{d+1}}{d+1} + 2(0.8n)^d + (0.2)(0.8)^d
          n^{d+1}.
    \end{align*}
    Consequently,
    \begin{align*}
      \sum_{j=1}^n j^d - \sum_{c \in \lambda/\nu} h_c^d
          &> \frac{n^{d+1}}{d+1} - \frac{(0.8n)^{d+1}}{d+1}
               - 2(0.8n)^d - (0.2)(0.8)^d n^{d+1} \\
          &= \left(\frac{1}{d+1} (1 - (0.8)^{d+1}) - 0.2 (0.8)^d\right) n^{d+1}
               - 2(0.8)^d n^d.
    \end{align*}
    It is easy to check that the coefficient on $n^{d+1}$ is bounded below by
    $\frac{1}{26(d+1)}$ for all positive integers $d$. The result follows.
  \end{proof}
\end{Lemma}

\begin{Definition}
  Given any partition $\lambda/\nu \vdash n$,
  let the \textit{aft} of $\lambda/\nu$ be  the statistic
    \[ \aft(\lambda/\nu)
        \coloneqq n - \max_{c \in \lambda/\nu} \{\arm(c), \leg(c)\} \]
  where $\arm(c)$ is the number of cells in the same row as $c$ to
  the right of $c$, including $c$ itself, and $\leg(c)$ is the number of
  cells in the same column as $c$ below $c$, including $c$. When
  $\nu = \varnothing$, we have
  $\aft(\lambda) = n - \max\{\lambda_1, \lambda_1'\}$
  as above. When $\lambda/\nu = \underline{\lambda}$,
  we have $\aft(\underline{\lambda}) = n -
  \max_i\{\lambda_1^{(i)}, {\lambda^{(i)}}'_1\}$.
  Note that $h_c = \arm(c) + \leg(c) - 1$.
\end{Definition}

\begin{Lemma}\label{lem:large_hook_bounds}
  Let $\lambda/\nu \vdash n$ such that $\max_{c \in \lambda/\nu} h_c
  \geq 0.8n$, and let $d$ be any positive integer. Furthermore, suppose
  $n \geq 10$. Then,
  \begin{align}
     \aft(\lambda/\nu) \frac{\lfloor 0.1n\rfloor^d}{d}  \hspace{.2cm} \leq  \hspace{.2cm}
    \sum_{j=1}^n j^d - \sum_{c \in \lambda/\nu} h_c^d  \hspace{.2cm}
    \leq  \hspace{.2cm} 2\aft(\lambda/\nu) \left(n^d + dn^{d-1}\right).
  \end{align}

  \begin{proof}
    The result holds trivially if $\aft(\lambda/\nu) = 0$ since
    in that case $\lambda/\nu$ is a single row or column, so assume
    $\aft(\lambda/\nu)>0$. Let $m \in \lambda/\nu$ have
    $h_m \geq 0.8n$, where we may assume $m$ is the
    first cell in its row and column. For convenience, we may further
    assume by symmetry that $\arm(m) \geq \leg(m)$.
    Since $h_m \geq 0.8n$, it also follows that $\aft(\lambda/\nu)
    = n - \arm(m)$.

    Now let $R$ be the set of cells in the row of $m$, not including
    $m$ itself, which are the only cells of $\lambda/\nu$ in their
    columns. Since $\lambda/\nu$ is a skew partition, $R$ is
    connected.  We claim that $\#R \geq 0.1n$. To prove the claim, we
    first observe that the hypothesis $h_m \geq 0.8n$ implies there
    are at most $n-h_m \leq 0.2n$ cells of $\lambda/\nu$ which could
    possibly be in the columns of the cells of the row of $m$ not
    including $m$.  Since $\arm(m) \geq \leg(m)$ and
    $\arm(m) + \leg(m) - 1 = h_m \geq 0.8n$, we have
    $\arm(m) \geq 0.4n$. Hence no more than $0.2n$ of the $0.4n-1$
    cells in the row of $m$ not including $m$ can be excluded from
    $R$, so $\#R \geq 0.4n-1-0.2n \geq 0.1n$ for $n \geq 10$.

    Construct $T \in \RSYT(\lambda/\nu)$ iteratively as follows; see
    \Cref{fig:rsyt_filling.2} for an example.
    At each step of the iteration, we will first increment all existing labels
    by $1$ and then label a new outer cell with $1$. Begin by adding the cells of
    the row of $m$ from left to right until the last cell of $R$ has been
    added. Now add the remaining cells of $\lambda/\nu$ row by row
    starting at the topmost row and going from left to right. It is easy
    to see that the result respects the decreasing row and column conditions, so
    $T \in \RSYT(\lambda/\nu)$.

    \begin{figure}[ht]
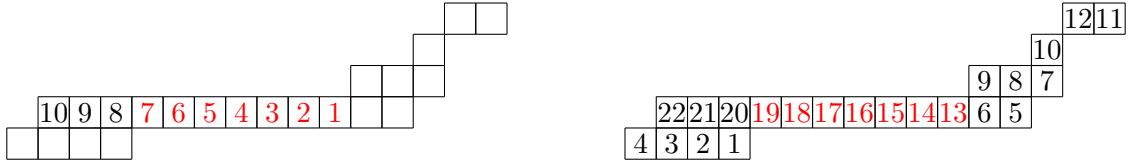

      \centering
      \begin{subfigure}[t]{0.49\textwidth}
        \centering
        \[
          \tableau{  &  &  &  &  &  &  &  &  &  &  &  &  &  &\ &\ \\
                          &  &  &  &  &  &  &  &  &  &  &  &  &\ &  &  \\
                          &  &  &  &  &  &  &  &  &  &  &\ &\ &\ &  &  \\
                          &10&9&8&\red 7&\red 6&\red 5&\red 4&\red 3&\red 2&\red 1&\ &\ \\
                       \ &\ &\ &\ }
        \]
      \end{subfigure}
      \begin{subfigure}[t]{0.49\textwidth}
        \centering
        \[
          \tableau{  &  &  &  &  &  &  &  &  &  &  &  &  &  &12&11\\
                          &  &  &  &  &  &  &  &  &  &  &  &  &10&  &  \\
                          &  &  &  &  &  &  &  &  &  &  &9&8&7&  &  \\
                          &22&21&20&\red{19}&\red{18}&\red{17}&\red{16}&\red{15}&\red{14}&\red{13}&6&5\\
                       4&3&2&1}
        \]
      \end{subfigure}
      \caption{On the left, the partially constructed $T \in \RSYT(\lambda/\nu)$
        after all the cells of $R$ (in red) have been filled. On the right,
        the final $T \in \RSYT(\lambda/\nu)$. Here $\aft(\lambda/\nu) = 10$.}
      \label{fig:rsyt_filling.2}
    \end{figure}

    By \Cref{lem:rsyt}, we have inequalities $T_c \geq h_c$.
    At every step of the iteration, a labeled cell
    has $T_c$ increase by $1$, while $h_c$ increases by $1$ if
    and only if the newly  labeled cell is in the hook of $c$. That is,
    for the final filling $T$, $T_c - h_c$ counts the number of
    times after cell $c$ was filled that the new cell was not in
    the same row or column as $c$. For each  $c \in R$,
    it follows that
    $T_c - h_c = n - \arm(m) = \aft(\lambda/\nu)$.

    For the lower bound, we now find
    \begin{align*}
      \sum_{k=1}^n k^d - \sum_{c \in \lambda/\nu} h_c^d
        &= \sum_{c \in R} (T_c - h_c){\mathbf h}_{d-1}(T_c, h_c) \\
        &= \sum_{c \in R} \aft(\lambda/\nu){\mathbf h}_{d-1}
             (h_c + \aft(\lambda/\nu), h_c) \\
        &\geq \sum_{k=1}^{\lfloor 0.1n\rfloor} \aft(\lambda/\nu)
             {\mathbf h}_{d-1}(k + \aft(\lambda/\nu), k) \\
        &\geq \aft(\lambda/\nu) \sum_{k=1}^{\lfloor 0.1n\rfloor} k^{d-1} \\
        &\geq \aft(\lambda/\nu) \frac{\lfloor 0.1n\rfloor^d}{d},
    \end{align*}
        where the first inequality uses the fact that
    $\{h_c : c \in R\}$ has pointwise lower bounds of
    $\{1, 2, \ldots, \#R\}$ and the last inequality uses
    \eqref{eq:sum.bounds}.

    For the upper bound, we construct a new $T \in \RSYT(\lambda/\nu)$
    as follows; see \Cref{fig:rsyt_filling} for an example. First, for each cell $c$
    in the row of $m$ taken
    from left to right, add the topmost cell in the column of $c$. Now add the
    remaining cells of $\lambda/\nu$ exactly as before. Again consider the
    final differences $T_c - h_c$. For cells added in the second stage,
    $T_c - h_c$ could increase no more than $n-\arm(m) = \aft(\lambda/\nu)$
    times, so $T_c - h_c \leq \aft(\lambda/\nu)$ for such $c$. For cells
    added in the first stage, we claim that $T_c - h_c \leq 2\aft(\lambda/\nu)$.
    For the claim, it suffices to show that after the first stage, for cells
    added in the first stage, $T_c - h_c \leq \aft(\lambda/\nu)$. During
    the first stage, the differences $T_c - h_c$ are zero while
    cells of row $m$ are being added. Afterwards during the first phase,
    cells not in row $m$ are added, of which there are no more than
    $n-\arm(m) = \aft(\lambda/\nu)$, so the differences $T_c - h_c$
    can increase no more than $\aft(\lambda/\nu)$ many times during the
    first phase, completing the claim.

    \begin{figure}[ht]
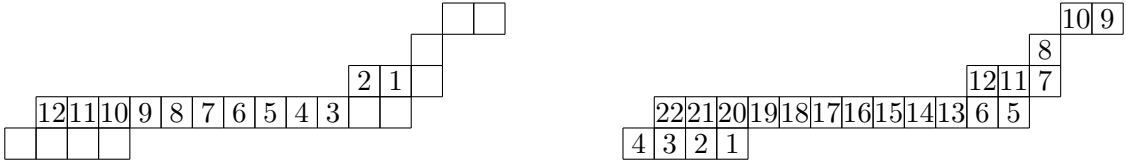

      \centering
      \begin{subfigure}[t]{0.49\textwidth}
        \centering
        \[
          \tableau{  &  &  &  &  &  &  &  &  &  &  &  &  &  &\ &\ \\
                          &  &  &  &  &  &  &  &  &  &  &  &  &\ &  &  \\
                          &  &  &  &  &  &  &  &  &  &  &2&1&\ &  &  \\
                          &12&11&10&9&8&7&6&5&4&3&\ &\ \\
                       \ &\ &\ &\ }
        \]
      \end{subfigure}
      \begin{subfigure}[t]{0.49\textwidth}
        \centering
        \[
          \tableau{  &  &  &  &  &  &  &  &  &  &  &  &  &  &10&9\\
                          &  &  &  &  &  &  &  &  &  &  &  &  &8&  &  \\
                          &  &  &  &  &  &  &  &  &  &  &12&11&7&  &  \\
                          &22&21&20&19&18&17&16&15&14&13&6&5\\
                        4&3&2&1}
        \]
      \end{subfigure}
      \caption{On the left, the second partially constructed $T \in \RSYT(\lambda/\nu)$
        after the first $\arm(m)$ cells have been filled. On the right,
        the final $T \in \RSYT(\lambda/\nu)$.}
      \label{fig:rsyt_filling}
    \end{figure}

    Having established that $T_c - h_c \leq 2\aft(\lambda/\nu)$,
    we now find by \eqref{eq:rsyt2} and \eqref{eq:sum.bounds},
    \begin{align*}
      \sum_{k=1}^n k^d - \sum_{c \in \lambda/\nu} h_c^d
        &= \sum_{c \in \lambda/\nu} (T_c - h_c) {\mathbf h}_{d-1}(T_c, h_c) \\
        &\leq \sum_{c \in \lambda/\nu} 2\aft(\lambda/\nu) {\mathbf h}_{d-1}(T_c, T_c) \\
        &= 2\aft(\lambda/\nu) \sum_{j=1}^n dj^{d-1} \\
        &< 2\aft(\lambda/\nu) \left(n^d + dn^{d-1}\right).
    \end{align*}
  \end{proof}
\end{Lemma}

\begin{Corollary}\label{cor:aft_theta}
  For fixed $d \in \bZ_{\geq 1}$, uniformly for all skew shapes $\lambda/\nu$,
    \begin{equation} \sum_{k=1}^{|\lambda/\nu|} k^d - \sum_{c \in \lambda/\nu} h_c^d
      = \Theta(\aft(\lambda/\nu) \cdot |\lambda/\nu|^d). \end{equation}

    % That is, there are constants $c_1, c_2 > 0$ such that
%    \[ c_1 \aft(\lambda/\nu) n^d \leq
%        \sum_{k=1}^n k^d - \sum_{c \in \lambda/\nu} h_c^d
%        \leq c_2 \aft(\lambda/\nu) n^d. \]

  \begin{proof}
    Let $n = |\lambda/\nu|$.
    When $\max_{c \in \lambda/\nu} h_c \geq 0.8n$, the result
    follows from \Cref{lem:large_hook_bounds}. On the other hand, when
    $\max_{c \in \lambda/\nu} h_c < 0.8n$, then
    $n \geq \aft(\lambda/\nu) \geq 0.2n$, and the result follows
    from \Cref{lem:small_hook_bounds}.
  \end{proof}
\end{Corollary}

\begin{Corollary}\label{cor:bound.cumulants}
  Fix $d$ to be an even positive integer. Uniformly for
  all block diagonal skew shapes $\underline{\lambda}$,
  the absolute value of the normalized cumulant
  $|{\kappa_d^{\underline{\lambda}}}^*|$ of
  $\cX_{\underline{\lambda}}[\maj]$ is $\Theta(\aft(\underline{\lambda})^{1-d/2})$.

  \begin{proof}
    For $d$ even, by \eqref{eq:cumulant_diag} and \Cref{cor:aft_theta},
    we have
      \[ |\kappa_d^{\underline{\lambda}}|
         = \Theta(\aft(\underline{\lambda}) n^d), \]
    where $n = |\underline{\lambda}|$.
    Consequently by the homogeneity of cumulants, we have
      \[ |{\kappa_d^{\underline{\lambda}}}^*|
          = \left|\frac{\kappa_d^{\underline{\lambda}}}
             {(\kappa_2^{\underline{\lambda}})^{d/2}}\right|
          = \Theta\left(\frac{\aft(\underline{\lambda})n^d}
             {\aft(\underline{\lambda})^{d/2} n^d}\right)
          = \Theta(\aft(\underline{\lambda})^{1-d/2}). \]
  \end{proof}
\end{Corollary}

We now state and prove the generalization of \Cref{thm:an} for the block
diagonal skew shapes $\underline{\lambda}$ from
\Cref{ssec:words_tableaux}.

\begin{Theorem}\label{thm:an_diag}
  Suppose $\underline{\lambda}^{(1)}, \underline{\lambda}^{(2)}, \ldots$
  is a sequence of block diagonal skew partitions, and let
  $\cX_N \coloneqq \cX_{\underline{\lambda}^{(N)}}[\maj]$ be the corresponding
  random variables for the $\maj$ statistic. Then, the sequence $\cX_1,
  \cX_2, \ldots$ is asymptotically normal if and only if
  $\aft(\underline{\lambda}^{(N)}) \to \infty$ as $N \to \infty$.

  \begin{proof}
    If $\aft(\underline{\lambda}^{(N)}) \to \infty$,
    the result follows immediately from \Cref{cor:cumulants},
    \Cref{cor:bound.cumulants}, and the fact that the
    odd cumulants vanish. On the other hand, if
    $\aft(\underline{\lambda}^{(N)}) \not\to \infty$,
    in the next  section we will show that $\cX_1^*, \cX_2^*, \ldots$
    has a subsequence which converges to either a discrete
    or uniform-sum distribution, which in either case is non-normal.
  \end{proof}
\end{Theorem}

\begin{Remark}
  Using work of Hwang--Zacharovas \cite[Thm.~1.1]{MR3346464},
  considering just the $d=4$ case is sufficient to prove both
  directions of \Cref{thm:an_diag}. However, the estimates we've given
  for $\kappa_d^{\underline{\lambda}}$ are strong enough to bound
  all the normalized cumulants simultaneously, and restricting to $d=4$
  (or even $d=2$) does not simplify the argument.
\end{Remark}

%% todo:  add comment from Janson converastation:   "Yes, I forgot
%% your blanket assumption. It does indeed take care of most
%% problems. But in order to conclude convergence of moments when you
%% have convergence in distribution, you would need to assume some
%% uniformity in your exponential tail condition (after
%% normalization), so that the d:th moment is uniformly bounded."

\section{Uniform sum limits for $\maj$ on $\SYT(\protect\underline{\lambda})$}\label{sec:SYT_diag_IH}

The estimates from \Cref{sec:SYT_diag_an} apply when
$\aft \to \infty$. We next give an analogous estimate handling the
case when $\aft$ is bounded, resulting in \Cref{thm:asum_diag}.
We may then deduce \Cref{thm:all_limits} from the introduction
and its generalization to block diagonal skew shapes,
\Cref{thm:all_limits_diag}.
Recall from \Cref{sec:intro} and \Cref{ex:Ucont}
that $\cIH_M$ is the
\textit{Irwin--Hall} distribution obtained by adding $M$ i.i.d.~$\cU[0, 1]$
random variables.

\begin{Lemma}\label{lem:fixed_aft_bound}
  Suppose $\lambda^{(N)}/\nu^{(N)} \vdash n_N$ is a sequence of skew
  partitions such that $\lim_{N \to \infty} n_N = \infty$ and
    \begin{equation} \lim_{N \to \infty} \aft(\lambda^{(N)}/\nu^{(N)}) =
        M \in \bZ_{\geq 0}. \end{equation}
  Then for each fixed $d \in \bZ_{\geq 1}$, we have
    \begin{equation} \lim_{N \to \infty} \frac{\sum_{k=1}^{n_N} k^d
        - \sum_{c \in \lambda^{(N)}/\nu^{(N)}} h_c^d}
        {M n_N^d} = 1. \end{equation}

  \begin{proof}
    Take $N$ large enough so that $\aft(\lambda^{(N)}/\nu^{(N)}) = M$
    and $n_N \gg M$. Let $m \in \lambda^{(N)}/\nu^{(N)}$ be such that
    $\aft(\lambda^{(N)}/\nu^{(N)}) = M = n_N - \arm(m)$ so $m$ is the
    first cell in its row and column, as in the proof of
    \Cref{lem:large_hook_bounds}.  Consider three regions of
    $\lambda^{(N)}/\nu^{(N)}$:
    \begin{enumerate}[(i)]
      \item The rightmost $\arm(m) - M = n_N-2M$ cells in the row of $m$.
      \item The remaining leftmost $M$ cells in the row of $m$.
      \item The remaining $M$ cells in $\lambda^{(N)}/\nu^{(N)}$.
    \end{enumerate}
    Construct $T \in \RSYT(\lambda^{(N)}/\nu^{(N)})$ iteratively
    as in the proof of \Cref{lem:large_hook_bounds} as follows.
    First add cells in region (iii) row by row starting at the topmost row proceeding
    from left to right, stopping just before inserting the row of $m$.
    Next add the cells from region (ii) from left to right.
    Now add the remaining cells in region (iii) row by row starting at the row immediately
    below the row of $m$ proceeding from left to right. Finally
    insert the cells from region (i) from left to right. It is easy
    to see that the cells in region (i) are the lowest cells in their
    column, from which it follows that $T$ indeed satisfies the
    column and row decreasing conditions.

    We now consider the contributions of regions (i)-(iii) to the
    quotient
    $$
    \frac{\sum_{k=1}^{n_N} k^d
        - \sum_{c \in \lambda^{(N)}/\nu^{(N)}} h_c^d}
        {M n_N^d}.
    $$
    Recall that $T_c - h_c$ can be interpreted as the
    number of times a cell inserted after cell $c$ was not inserted
    in the same hook as $c$. It follows that $T_c - h_c = 0$
    for region (i), leaving only contributions from the $2M$ cells
    in regions (ii) and (iii), a bounded sum. For region (ii), we have
    $T_c-h_c\leq M$, so that
      \[ T_c^d-h_c^d=(T_c - h_c){\mathbf h}_{d-1}(T_c, h_c) \leq (2M)dn_N^{d-1}. \]
    Dividing by $Mn_N^d$, cells in region (ii) contribute $0$ to the
    sum in the limit. Finally, for region (iii), we find
    $1 \leq h_c \leq M+1$ and $n_N - 2M+1 \leq T_c \leq n_N$,
    so that for each of the $M$ cells $c$ in region (iii),
      \[ (n_N - 2M+1)^d - (M+1)^d \leq T_c^d - h_c^d
          \leq n_N^d - 1^d. \]
    Dividing by $n_N^d$, both bounds are asymptotic
    to $1$ as $n_N \to \infty$. Adding up all $M$ such contributions,
    the result follows.
  \end{proof}
\end{Lemma}

\begin{Theorem}\label{thm:asum_diag}
  Suppose that $\underline{\lambda}^{(1)}, \underline{\lambda}^{(2)}, \ldots$ is
  a sequence of block diagonal skew partitions such that
  $\lim_{N \to \infty} |\underline{\lambda}^{(N)}| = \infty$ and
  $\aft(\underline{\lambda}^{(N)}) = M$ is constant.
  Let $\cX_N \coloneqq \cX_{\underline{\lambda}^{(N)}}[\maj]$ be the corresponding
  random variable for the $\maj$ statistic. Then $\cX_1^*, \cX_2^*, \ldots$
  converges in distribution to $\cIH_M^*$.

  \begin{proof}
    Using \Cref{eq:cumulant_diag} and \Cref{lem:fixed_aft_bound},
    we have for $d \geq 2$ that
    \begin{align*}
      \lim_{N \to \infty} (\kappa_d^{\underline{\lambda}^{(N)}})^*
      &= \lim_{N \to \infty} \frac{\kappa_d^{\underline{\lambda}^{(N)}}}
         {(\kappa_d^{\underline{\lambda}})^{d/2}} \\
      &= \lim_{N \to \infty} \frac{(B_d/d) \left(\sum_{k=1}^{n_N} k^d
           - \sum_{c \in \underline{\lambda}^{(N)}} h_c^d\right)}
           {(B_2/2)^{d/2} \left(\sum_{k=1}^{n_N} k^2
           - \sum_{c \in \underline{\lambda}^{(N)}} h_c^2\right)^{d/2}} \\
      &= \lim_{N \to \infty} \frac{(B_d/d)}{(B_2/2)^{d/2}}
           \frac{Mn_N^d}{(Mn_N^2)^{d/2}} \\
      &= \frac{(MB_d/d)}{(MB_2/2)^{d/2}}.
    \end{align*}
    From \Cref{ex:Ucont} and the homogeneity and additivity properties
    of cumulants, we have
    \begin{align*}
      (\kappa_d^{\cIH_M})^*
        &= \frac{\kappa_d^{\cIH_M}}{(\kappa_2^{\cIH_M})^{d/2}} \\
        &= \frac{(MB_d/d)}{(MB_2/2)^{d/2}}.
    \end{align*}
    The result now follows from  \Cref{thm:moments} after converting
    moments to cumulants.
  \end{proof}
\end{Theorem}

\begin{Theorem}\label{thm:all_limits_diag}
  Let $\underline{\lambda}^{(1)}, \underline{\lambda}^{(2)}, \ldots$ be
  a sequence of block diagonal skew partitions.
  Then the sequence $(\cX_{\underline{\lambda}^{(N)}}[\maj]^*)$ converges in
  distribution if and only if
  \begin{enumerate}[(i)]
    \item $\aft(\underline{\lambda}^{(N)}) \to \infty$; or
    \item $|\underline{\lambda}^{(N)}| \to \infty$ and
      $\aft(\underline{\lambda}^{(N)}) \to M < \infty$; or
    \item the distribution of $\cX_{\underline{\lambda}^{(N)}}[\maj]$
    is eventually constant.
  \end{enumerate}
  The limit law is $\cN$ in case (i), $\cIH_M^*$ in case (ii),
  and discrete in case (iii).

  \begin{proof}
    The backwards direction follows from \Cref{thm:an_diag} and
    \Cref{thm:asum_diag}. In the forwards direction, let
    $\underline{\lambda}^{(N)}$ be such a sequence where
    $(\cX_{\underline{\lambda}^{(N)}}[\maj]^*)$ converges in
    distribution. If $|\underline{\lambda}^{(N)}|$ is bounded, then
    there are only finitely many distinct $\underline{\lambda}^{(N)}$,
    forcing case (iii).  If $|\underline{\lambda}^{(N)}|$ is
    unbounded, then we have subsequences satisfying either (i) or (ii)
    since the sequence converges in distribution, which from
    \Cref{thm:an_diag} and \Cref{thm:asum_diag} gives convergence in
    distribution to $\cN$ or $\cIH_M^*$, which are continuous,
    distinct distributions. The result follows.
  \end{proof}
\end{Theorem}

\begin{figure} [ht]
  \begin{center}
    \includegraphics[height=2.5in, trim=0 50 0 0, clip]{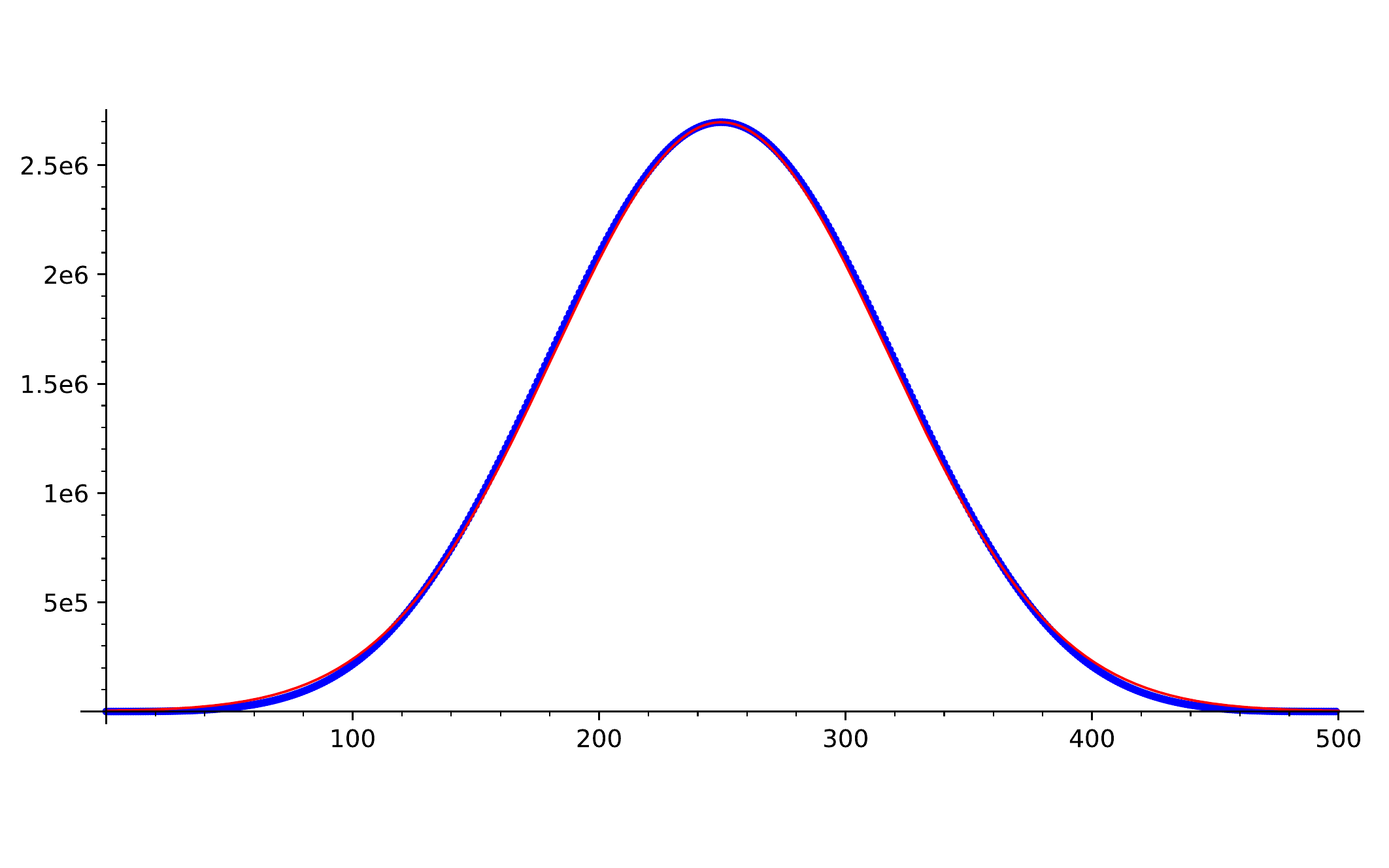}
    \caption{Coefficients of $\SYT(\lambda)^{\maj}(q)$ for
      $\lambda = (100,3,2)$ where $\aft(\lambda) = 5$ plotted in blue
      along with the corresponding normal distribution with the same
      mean and variance plotted in red.  The difference is mostly in
      the tails.}
    \label{fig:100-3-2}
  \end{center}
\end{figure}

From the Central Limit Theorem, we know the Irwin--Hall distribution
$\cIH_M^*$ for $M$ large closely resembles a normal distribution, so
it will be quite rare for a plot of the coefficients of
$\SYT(\lambda)^{\maj}(q)$ to look anything but normal.  Since Irwin--Hall
distributions are finitely supported, the difference between the two
distributions is mainly in the tails. We note that
even for $M=5$, there is a close resemblance. See the plot in
\Cref{fig:100-3-2}.

\section{Discrete distributions for $\maj$ on $\SYT(\lambda)$}\label{sec:SYT_discrete}

We conclude by analyzing more carefully the
discrete case of the limit law classification for
$\maj$ on $\SYT(\lambda)$, \Cref{thm:all_limits}. The result
is \Cref{thm:discrete_distributions}, which lists several families of
pairs of shapes $\lambda$ and $\nu$ of differing sizes
for which we nonetheless have $\#\SYT(\lambda) = \#\SYT(\nu)$.

A well-known corollary of \eqref{eq:SYT_maj_hook} is that for
partitions $\lambda$ and $\nu$ of $n$, $\maj$ is equidistributed on
$\SYT(\lambda)$ and $\SYT(\nu)$ if and only if $\rank(\lambda) = \rank(\nu)$
and the multisets $\{h_c : c \in \lambda\}$ and $\{h_d : d \in \nu\}$
are equal. These hook multisets do not entirely characterize the
partition---see \cite{MR0463279}.  The following theorem gives a
similar result even if we consider the corresponding standardized
random variables $\cX_\lambda[\maj]$ and $\cX_\nu[\maj]$.

\begin{Theorem}\label{thm:discrete_distributions}
  Let $\lambda$ and $\nu$ be partitions. Then $\cX_\lambda[\maj]^*$
  and $\cX_\nu[\maj]^*$ have the same distribution if and only if
  \begin{enumerate}[(i)]
    \item the multisets of hook lengths $\{h_c : c \in \lambda\}$ and $\{h_d : d \in \nu\}$
      are equal; or
    \item the multisets $\{h_c : c \in \lambda\}$ and
      $\{|\lambda|\} \sqcup \{h_d : d \in \nu\}$ are equal; or
    \item $\lambda$ and $\nu$ are each either a single row or column; or
    \item $\lambda, \nu \in \{(2, 1), (2, 2)\}$.
  \end{enumerate}
Moreover, case (ii) occurs if and only if, up to transposing,
  \begin{enumerate}[(a)]
    \item $\lambda = (n)$ and $\nu = (n-1)$ for $n \geq 2$; or
    \item $\lambda = (r+1, 1^{2r+2})$ and $\nu = (2^{r+1}, 1^r)$ for
      $r \geq 1$; or
    \item $\lambda = (s, 1^{s+2})$ and $\nu = (s, s, 1)$ for $s \geq 4$; or
    \item $\lambda = (3, 1^5)$ and $\nu = (3^2, 1)$,
      or $\lambda = (4, 1^6)$ and $\nu = (3^3, 1)$.
  \end{enumerate}

  \begin{proof}
    Let $n \coloneqq |\lambda|$ and $m \coloneqq |\nu|$.  Let
    $f^\lambda(q) = \frac{[n]_q!}{\prod_{c \in \lambda} [h_c]}$, which
    is a polynomial by \eqref{eq:SYT_maj_hook} with constant coefficient $1$.
    Let $f^\lambda=f^\lambda(1)=|\SYT(\lambda)|$. Let $f^{\nu}$ and
    $f^{\nu}(q)$ be defined similarly.

    In the backwards direction, if (i) holds, then $n=m$,
    both variances agree  by
    \Cref{thm:SYT_moments}, and $f^\lambda(q) = f^\nu(q)$, so
    $\cX_\lambda[\maj]^*$ and $\cX_\nu[\maj]^*$ have the same
    distribution.  Similarly if (ii) holds $f^\lambda(q) = f^\nu(q)$,
    both variances agree, and $\cX_\lambda[\maj]^*$ and $\cX_\nu[\maj]^*$ have
    the same distribution again. Condition (iii) holds if and only if
    the distributions are concentrated at a single point. For (iv), we
    have $f^{(2, 1)}(q) = 1+q$ and $f^{(2, 2)}(q) = 1+q^2$, so the
    normalized distributions are clearly equal.

    In the forwards direction, suppose $\cX_\lambda[\maj]^*$ and
    $\cX_\nu[\maj]^*$ have the same distribution.  Since $f^\lambda(q)$
    has constant coefficient $1$, $\cX_\lambda[\maj]$ is concentrated at
    a single point if and only if $f^\lambda=1$, which occurs if and
    only if $\lambda$ is a single row or column which is covered by
    case (iii).  It is easy to see that $f^\lambda = 2$ if and only if
    $\lambda \in \{(2, 1), (2, 2)\}$ which is covered by case (iv).

    Assume $f^\lambda, f^\nu > 2$. By \cite[Thm.~1.1]{BKS-zeroes},
    it follows that $f^\lambda(q)$ and $f^\nu(q)$ each have two adjacent
    non-zero coefficients.
    Since $f^\lambda(q)$ and $f^\nu(q)$ each have constant term 1 and
    two adjacent non-zero coefficients, then it follows from the
    assumption $\cX_\lambda[\maj]^*$ and $\cX_\nu[\maj]^*$ have the same
    distribution that
    \begin{equation}\label{eq:dhp}
      f^{\lambda}(q)=\frac{[n]_q!}{\prod_{c \in \lambda} [h_c]_q}
        = \frac{[m]_q!}{\prod_{d \in \nu} [h_d]_q} = f^{\nu}(q).
    \end{equation}

    Without loss of generality, we can assume $n \geq m$. If $n=m$, we
    have $\prod_{c \in \lambda} [h_c]_q = \prod_{d \in \nu} [h_d]_q$,
    from which it follows that the multisets of hook lengths are equal
    by considering multiplicities of zeros at all primitive roots of
    unity as in case (i).

    From here on, assume $n>m$.  The multiplicity of a zero of a
    primitive $n$th root of unity in \eqref{eq:dhp} is $0$ on the
    right, so from the left $\lambda$ must have a hook of length $n$
    so it itself is a hook shape partition. Since $\lambda$ is not a
    single row or column by the assumption $f^\lambda>2$, we know
    $\lambda$ does not have a cell with hook length $n-1$.
    Consequently,  the multiplicity of a zero at a primitive $(n-1)$th
    root of unity in \eqref{eq:dhp} is $1$ on the left, forcing
    $m=n-1$ on the right. Thus \eqref{eq:dhp} becomes
    \begin{equation}\label{eq:caseii}
    [m+1]_q \prod_{d \in \nu} [h_d]_q = \prod_{c \in \lambda}
    [h_c]_q,
    \end{equation}
    and as before the multiset condition (ii) must hold.
    This completes the proof of the first statement in the theorem.

    For the second statement, suppose (ii) holds, so the multisets
    $\{h_c : c \in \lambda\}$ and
    $\{|\lambda|\} \sqcup \{h_d : d \in \nu\}$ are equal.  Then,
    $m=n-1$ and $\lambda$ has a cell with hook length $|\lambda|$, so
    $\lambda$ is a hook shape partition $(n-k,1^{k})$ for some
    $0 \leq k \leq n$, and
    \begin{equation}\label{eq:dhp2}
      \{h_d : d \in \nu\} = [m-k] \sqcup [k].
    \end{equation}
    By transposing if necessary, we may assume $k \geq m-k$ is
    the maximum hook length in $\nu$.
    If $\lambda$ has one cell with hook length $1$, then (a) holds.
    Otherwise, both $\lambda$ and $\nu$ have precisely two cells with
    hook length $1$, so $\nu$ is the union of two rectangles and not
    itself a rectangle. If $\nu$ were a hook, then it would have a hook
    length equal to $m$ which would imply $\lambda$ has a cell of hook
    length $m=n-1$ contradicting the fact that $\lambda$ has two outer
    corners. Thus $\nu$ is not itself a hook.

    Transposing $\nu$ if necessary, we can assume its first two rows are equal,
    say $\nu_1 = \nu_2 = s$. If $\nu_1' = \nu_2'$, one may check
    that the cell furthest from the origin in the intersection of the two
    rectangles forming $\nu$ would be the only cell of its hook length, and
    that moreover its two neighbors in the intersection would each have one
    larger hook length, contrary to \eqref{eq:dhp2}. It follows that
    $\nu = (s^t, 1^r)$ where $r \geq 1$, $s \geq 2$, and $t \geq 2$.
    We now have several cases.
    \begin{itemize}
      \item If $s = 2$, the hook lengths of $\nu$ are
        $\{1, \ldots, r, r+2, \ldots, r+t+1, 1, \ldots, t\}$.
        The ``gap'' between $r$ and $r+2$ together with
        \eqref{eq:dhp2} forces $t=r+1$, so that
        $\nu = (2^{r+1}, 1^r)$ with $r \geq 1$.
        Here $k=r+t+1=2r+2$, resulting in case (b).
      \item If $s \geq 3$, the last two columns of $\nu$ already contain two cells
        with hook length $2$. If $r>1$, the first column would also have a cell
        with hook length $2$, contradicting \eqref{eq:dhp2}, so $r=1$.
        \begin{itemize}
          \item If $s = 3$, the hook lengths of $\nu$ are
            $\{1, \ldots, t, 2, \ldots, t+1, 1, 4, 5, \ldots, t+3\}$. Because of the
            ``gap'' between $t+1$ and $t+3$, this is of the form in
            \eqref{eq:dhp2} if and only if $t=2$ or $t=3$, resulting in case (d).
          \item Suppose $s > 3$. If $t \geq 3$, then the final three columns of $\nu$
            contain three cells with hook length $3$, contradicting \eqref{eq:dhp2}, so $t=2$.
            The hook lengths of $\nu$ are then
            $\{1, 1, 2, \ldots, s-1, s+1, 2, 3, \ldots, s, s+2\}$, which is already of the form
            \eqref{eq:dhp2}, resulting in case (c).
        \end{itemize}
    \end{itemize}
    The reverse implications from (a)-(d) to (ii) were verified in the course of the
    above argument.
  \end{proof}
\end{Theorem}

\begin{Remark}
  The proof of \Cref{thm:discrete_distributions} applies more generally to arbitrary
  scaling factors and translations of the distributions of $\cX_\lambda[\maj]$ and
  $X_\nu[\maj]$, and not just those coming from means and variances.
\end{Remark}

\section{Future work}\label{sec:future}

We conjecture that almost all of the polynomials of the form
$\SYT(\lambda)^{\maj}(q)$ are unimodal and log-concave.  In this section, we
discuss the deviations of each of these properties.
In the rare cases where unimodality or
log-concavity fails, it only seems to happen at the very beginning
and end of the sequence of coefficients or near the middle
coefficient.

Recall that a polynomial $P(q) = \sum_{i=0}^n c_i q^i$ is
\emph{unimodal} if
\[ c_0\leq c_1 \leq \cdots \leq c_j \geq c_{j+1} \geq \cdots \geq
  c_n \] for some $j$, and $P(q)$ is \emph{log-concave} if
$c_i^2 \geq c_{i-1}c_{i+1}$ for all integers $0<i<n$.  A polynomial
with nonnegative coefficients which is log-concave and has no internal
zero coefficients is necessarily unimodal \cite{Stanley.1989}.
By \cite{BKS-zeroes}, we know exactly where internal zeros occur so
log-concavity would imply unimodality in these cases.

We say $P(q)$ is \emph{nearly unimodal} if instead
  \[ c_0\leq c_1 \leq \cdots \leq c_j, c_{j+1} = c_j - 1 < c_{j+2}
    \leq \cdots \leq c_{\lfloor \frac{n}{2} \rfloor}
  \]
  for some $j$ and $P(q)$ has symmetric coefficients.  Also, a
  symmetric polynomial $P(q)$ is \emph{nearly log-concave} if
  $c_i^2 \geq c_{i-1}c_{i+1}$ for all
  $1<i<\lfloor \frac{n}{2} \rfloor$.

  \begin{Conjecture}\label{conj:unimodality}
  The polynomial $\SYT(\lambda)^{\maj}(q)$ is unimodal if $\lambda$ has at least
  $4$ corners. If $\lambda$ has $3$ corners or fewer, then $\SYT(\lambda)^{\maj}(q)$
  is unimodal except when $\lambda$ or $\lambda'$ is among the following partitions:
  \begin{enumerate}
  \item Any partition of rectangle shape that has more than one row
    and column.

  \item Any partition of the form $(k,2)$ with $k\geq 4$ and $k$
    even.

  \item Any partition of the form $(k,4)$ with $k\geq 6$ and $k$
    even.

  \item Any partition of the form $(k,2,1,1)$ with $k\geq 2$ and $k$
    even.

  \item Any partition of the form $(k,2,2)$ with $k\geq 6$.

  \item Any partition on the list of 40 special exceptions:
    \begin{gather*}
      (3, 3, 2), (4, 2, 2), (4, 4, 2), (4, 4, 1, 1), (5, 3, 3),
        (7, 5), (6, 2, 1, 1, 1, 1), \\
      (5, 5, 2), (5, 5, 1, 1), (5, 3, 2, 2), (4, 4, 3, 1), (4, 4, 2, 2),
      (7, 3, 3), (8, 6), (6, 6, 2), \\
      (6, 6, 1, 1), (5, 5, 2, 2), (5, 3, 3, 3), (4, 4, 4, 2),
        (11, 5), (10, 6), (9, 7), (7, 7, 2), \\
      (7, 7, 1, 1), (6, 6, 4), (6, 6, 1, 1, 1, 1), (6, 5, 5), (5, 5, 3, 3),
        (12, 6), (11, 7), (10, 8), \\
      (15, 5), (14, 6), (11, 9), (16, 6), (12, 10), (18, 6),
        (14, 10), (20, 6), (22, 6).
     \end{gather*}
  \end{enumerate}
\end{Conjecture}

\Cref{conj:unimodality} was checked for all partitions up to size
$n= 50$.  Each of the families $(k,2)$,\ $(k,4)$,\ or $(k,2,1,1)$ have
a relatively simple set of hook lengths so explicit formulas can be
derived for the coefficients of $\SYT(\lambda)^{\maj}(q)$.  We have
found explicit proofs of near unimodality for each of these cases.
They are related to known integer sequences \cite[A266755]{oeis} and
\cite[A008642]{oeis} with nice generating functions.  Furthermore,
these families are all nearly unimodal as well as 20 of the special
exceptions.  All rectangles with at least 2 rows and columns are
nearly unimodal for $30 \leq n \leq 100$.  The only deviation occurs at $i=1$ up
to symmetry.  We conjecture this trend
also continues, hence the claim that all coefficients in
$\SYT(\lambda)^{\maj}(q)$ are close to unimodal.  The family $(k,2,2)$
is a bit further from being unimodal. The proof of the following
result is omitted, but follows directly from a careful analysis of the
hook lengths.

\begin{Proposition}
  If $\lambda=(k,2,2)$ for any positive integer $k \geq 3$, then the maximal
  coefficient of $f^\lambda(q)$, say $c_j$, satisfies
  the equation $c_{j} = c_{j+1} + \mathrm{floor}(k/6) + I(4 = (k\ \mathrm{mod}\ 6))$ and
  $c_0 \leq c_1 \leq \cdots \leq c_{j}$ and $j+1$ is the median
  nonzero coefficient.  Here $I$ is an indicator
  function which is 1 if true and 0 if false.
\end{Proposition}

 \begin{Conjecture}\label{conj:nearly-unimodal}
   The polynomials $\SYT(\lambda)^{\maj}(q)$ are ``nearly unimodal but not
   unimodal'' for partitions $\lambda$ or $\lambda'$ in the following
   cases:
   \begin{enumerate}
          \item Any partition of rectangle shape that has more than one row
     and column with more than 30 cells.

   \item Any partition of the form $(k,2)$ with $k\geq 4$ and $k$
     even.

   \item Any partition of the form $(k,4)$ with $k\geq 6$ and $k$
     even.

   \item Any partition of the form $(k,2,1,1)$ with $k\geq 2$ and $k$
     even.
     \end{enumerate}
 \end{Conjecture}

\Cref{conj:nearly-unimodal} was checked for all paritions of size up to $n=100$.  It also
 holds for the following 14 special exceptions:
 $$
   (3, 3, 2), (4, 2, 2), (5, 3, 3), (7, 5), (6, 2, 1, 1, 1, 1), (5, 3, 2, 2), (4,
   4, 3, 1), (7, 3, 3), (5, 3, 3, 3),$$
   $$  (11, 5), (6, 6, 1, 1, 1, 1), (6, 5, 5), (15, 5), (22, 6).$$

   \bigskip

Log-concavity for the polynomials $\SYT_\lambda^{\maj}(q)$ appears to
be harder to characterize.  There are examples of partitions with even
5 corners which are not log-concave.  For example $f^\lambda(q)$ for
$\lambda = (9, 9, 7, 7, 5, 5, 3, 3, 2)$ is nearly log-concave but
$c_1^2=4^2 = 16 <17 = c_0c_2$.  The only deviation occurs at $i=1$ up
to symmetry.
Thus, we summarize what we have observed in the following conjecture.

\begin{Conjecture}\label{conj:log-concave}
  The polynomials $\SYT(\lambda)^{\maj}(q)$ are almost always log-concave for
  partitions $\lambda \vdash n$ for large $n$.
\end{Conjecture}

This conjecture is based on the fact that the normal distribution is
log-concave and the following evidence.  The approximate probability
that a uniformly chosen partition of $n$ has the log-concave property
$\bP(\mathrm{LC})$ and the corresponding probability for the nearly
log-concave property $\bP(\mathrm{NLC})$ is given in the following
table:

\begin{figure}[ht]
  \centering
  \begin{tabular}{c|c|c|c}
    \toprule
    n & 30 & 40 & 50 \\
    \midrule
    \midrule
    $\bP(\mathrm{LC})$ & 0.6734475 & 0.7876426 & 0.8753587 \\
    $\bP(\mathrm{NLC})$ & 0.8003212 & 0.9204832 & 0.9688140 \\
    \bottomrule
  \end{tabular}
  \caption{Data supporting \Cref{conj:log-concave}.}
  \label{fig:log-concave}
\end{figure}

By \Cref{thm:an} and the conjectured claim that the coefficients of
$\SYT(\lambda)^{\maj}(q)$ are unimodal or almost unimodal for large
$\lambda$, one might hope that we could approximate the number of
$T \in \SYT(\lambda)$ with $\maj(T)=k$ by the density function
$f(k; \kappa_1^\lambda, \kappa_2^\lambda)$ for the normal distribution
with mean $\kappa_1^\lambda$ and variance $\kappa_2^\lambda$.  We have
the following conjectured bounds on such an approximation.

\begin{Conjecture}\label{conj:local.limit}
  Let $\lambda \vdash n$ be any partition. Uniformly for all $n$, for all integers $k$, we have
  \begin{equation*}
    \left|\bP[X_\lambda[\maj] = k] - f(k; \kappa_1^\lambda,
      \kappa_2^\lambda)\right| = O\left(\frac{1}{\sigma_\lambda \aft(\lambda)}\right).
  \end{equation*}
\end{Conjecture}

The conjecture has been verified for $25 < n \leq 50$ and $\aft(\lambda) > 1$ with a
constant of $1/9$, which is tight up to reasonable limits on computation in the sense that if it
is changed to $1/10$ with the other constraints the same, it fails at $n=50$.

\begin{Conjecture}\label{conj:skew_an}
  Asymptotic normality for general skew shapes and not just block
  diagonal skew shapes holds if and only if
  $\aft(\lambda/\nu^{(N)}) \to \infty$ as $N \to \infty$,
  generalizing the result in \Cref{thm:an_diag}.
\end{Conjecture}

%% older version:
% \begin{Conjecture}\label{conj:skew_an}
%   The asymptotic normality result \Cref{thm:an_diag} holds
%   for general skew shapes and not just block diagonal skew shapes.
% \end{Conjecture}

The argument in \Cref{sec:SYT_diag_an} proves that the
``formal cumulants'' associated with
  \[ \frac{[n]_q!}{\prod_{c \in \lambda/\mu} [h_c]_q} \]
exhibit asymptotic normality when $\aft(\lambda/\mu) \to \infty$.
However, this is only the first term in the general
$q$-Naruse formula for $\SYT(\lambda/\mu)^{\maj}(q)$. One
approach to \Cref{conj:skew_an} would be to show the remaining
terms are ``appropriately negligible.''

\section*{Acknowledgments}\label{sec:ack}

We would like to thank Krzysztof Burdzy, Rodney Canfield, Persi
Diaconis, Sergey Fomin, Pavel Galashin, Svante Janson, William
McGovern, Andrew Ohana, Greta Panova, Mihael Perman, Martin Rai\v{c},
Richard Stanley, Sheila Sundaram, Vasu Tewari, Lauren Williams, and
Alex Woo for helpful discussions related to this work.

\bibliography{syt_maj_2019}{}

\begin{thebibliography}{{Maz}08}

\bibitem[AR01]{MR1841639}
Ron~M. Adin and Yuval Roichman.
\newblock Descent functions and random {Y}oung tableaux.
\newblock {\em Combin. Probab. Comput.}, 10(3):187--201, 2001.

\bibitem[AS18]{AHLBACH201837}
Connor Ahlbach and Joshua~P. Swanson.
\newblock Refined cyclic sieving on words for the major index statistic.
\newblock {\em European Journal of Combinatorics}, 73:37 -- 60, 2018.

\bibitem[Bil95]{MR1324786}
Patrick Billingsley.
\newblock {\em Probability and measure}.
\newblock Wiley Series in Probability and Mathematical Statistics. John Wiley
  \& Sons, Inc., New York, third edition, 1995.
\newblock A Wiley-Interscience Publication.

\bibitem[BKS18]{BKS-zeroes}
Sara~C. {Billey}, Matja{\v{z}} {Konvalinka}, and Joshua~P. {Swanson}.
\newblock {Tableaux posets and the fake degrees of coinvariant algebras}.
\newblock {\em Preprint arXiv:1809.07386}, Sep 2018.

\bibitem[B{\'o}n15]{MR3408702}
Mikl\'os B{\'o}na, editor.
\newblock {\em Handbook of enumerative combinatorics}.
\newblock Discrete Mathematics and its Applications (Boca Raton). CRC Press,
  Boca Raton, FL, 2015.

\bibitem[Car75]{carlitz.1975}
L.~Carlitz.
\newblock A combinatorial property of {$q$}-{E}ulerian numbers.
\newblock {\em Amer. Math. Monthly}, 82:51--54, 1975.

\bibitem[CF13]{MR3084430}
Thomas Church and Benson Farb.
\newblock Representation theory and homological stability.
\newblock {\em Adv. Math.}, 245:250--314, 2013.

\bibitem[CJZ11]{MR2794017}
E.~Rodney Canfield, Svante Janson, and Doron Zeilberger.
\newblock The {M}ahonian probability distribution on words is asymptotically
  normal.
\newblock {\em Adv. in Appl. Math.}, 46(1-4):109--124, 2011.

\bibitem[CJZ12]{MR2925927c}
E.~Rodney Canfield, Svante Janson, and Doron Zeilberger.
\newblock Corrigendum to ``{T}he {M}ahonian probability distribution on words
  is asymptotically normal'' [{A}dv. in {A}ppl. {M}ath. 46 (1--4) (2011)
  109--124] [mr2794017].
\newblock {\em Adv. in Appl. Math.}, 49(1):77, 2012.

\bibitem[CWW08]{Chen-Wang-Wang.2008}
William Y.~C. Chen, Carol~J. Wang, and Larry X.~W. Wang.
\newblock The limiting distribution of the coefficients of the {$q$}-{C}atalan
  numbers.
\newblock {\em Proc. Amer. Math. Soc.}, 136(11):3759--3767, 2008.

\bibitem[DB62]{MR0155371}
F.~N. David and D.~E. Barton.
\newblock {\em Combinatorial chance}.
\newblock Hafner Publishing Co., New York, 1962.

\bibitem[Dia88]{MR964069}
Persi Diaconis.
\newblock {\em Group representations in probability and statistics}, volume~11
  of {\em Institute of Mathematical Statistics Lecture Notes---Monograph
  Series}.
\newblock Institute of Mathematical Statistics, Hayward, CA, 1988.

\bibitem[EL41]{MR0004841}
Paul Erd\"os and Joseph Lehner.
\newblock The distribution of the number of summands in the partitions of a
  positive integer.
\newblock {\em Duke Math. J.}, 8:335--345, 1941.

\bibitem[ER15]{Ehrenborg.Readdy.2015}
Richard Ehrenborg and Margaret Readdy.
\newblock A poset view of the major index.
\newblock {\em Adv. in Appl. Math.}, 62:1--14, 2015.

\bibitem[Fel45]{MR0013252}
W.~Feller.
\newblock The fundamental limit theorems in probability.
\newblock {\em Bull. Amer. Math. Soc.}, 51:800--832, 1945.

\bibitem[FH85]{MR814413}
J.~F\"urlinger and J.~Hofbauer.
\newblock {$q$}-{C}atalan numbers.
\newblock {\em J. Combin. Theory Ser. A}, 40(2):248--264, 1985.

\bibitem[Foa68]{Foata}
D.~Foata.
\newblock {On the Netto inversion number of a sequence}.
\newblock {\em Proc.\ Amer.\ Math.\ Soc.}, 19:236--2479--1130, 1968.

\bibitem[FRT54]{Frame-Robinson-Thrall.1954}
J.~S. Frame, G.~de~B. Robinson, and R.~M. Thrall.
\newblock The hook graphs of the symmetric groups.
\newblock {\em Canadian J. Math.}, 6:316--324, 1954.

\bibitem[Ful98]{MR1652841}
Jason Fulman.
\newblock The distribution of descents in fixed conjugacy classes of the
  symmetric groups.
\newblock {\em J. Combin. Theory Ser. A}, 84(2):171--180, 1998.

\bibitem[GKP89]{GKP}
Ronald~L. Graham, Donald~E. Knuth, and Oren Patashnik.
\newblock {\em Concrete mathematics}.
\newblock Addison-Wesley Publishing Company, Advanced Book Program, Reading,
  MA, 1989.

\bibitem[Gon44]{gon44}
V.~L. Goncharov.
\newblock From the realm of combinatorics.
\newblock {\em Izv. Akad. Nauk SSSR}, 8:3--48, 1944.

\bibitem[Har67]{MR0211432}
L.~H. Harper.
\newblock Stirling behavior is asymptotically normal.
\newblock {\em Ann. Math. Statist.}, 38:410--414, 1967.

\bibitem[HC78]{MR0463279}
Joan~E. Herman and Fan R.~K. Chung.
\newblock Some results on hook lengths.
\newblock {\em Discrete Math.}, 20(1):33--40, 1977/78.

\bibitem[HHL05]{HHL2005}
J.~Haglund, M.~Haiman, and N.~Loehr.
\newblock A combinatorial formula for {M}acdonald polynomials.
\newblock {\em J. Amer. Math. Soc.}, 18(3):735--761, 2005.

\bibitem[HZ15]{MR3346464}
Hsien-Kuei Hwang and Vytas Zacharovas.
\newblock Limit distribution of the coefficients of polynomials with only unit
  roots.
\newblock {\em Random Structures Algorithms}, 46(4):707--738, 2015.

\bibitem[IM65]{Iw.M}
N.~Iwahori and H.~Matsumoto.
\newblock On some {B}ruhat decomposition and the structure of the {H}ecke rings
  of {$\mathfrak{p}$}-adic {C}hevalley groups.
\newblock {\em Inst. Hautes \'{E}tudes Sci. Publ. Math.}, 25:5--48, 1965.

\bibitem[KL18]{1803.10457}
Gene~B. {Kim} and Sangchul {Lee}.
\newblock {Central limit theorem for descents in conjugacy classes of
  \$S\_n\$}.
\newblock {\em Preprint arXiv:1803.10457}, Mar 2018.

\bibitem[KO17]{MR3551643}
Jang~Soo Kim and Suho Oh.
\newblock The {S}elberg integral and {Y}oung books.
\newblock {\em J. Combin. Theory Ser. A}, 145:1--24, 2017.

\bibitem[Mac13]{MacMahon.1913}
P.~A. MacMahon.
\newblock The indices of permutations and the derivation therefrom of functions
  of a single variable associated with the permutations of any assemblage of
  objects.
\newblock {\em Amer. J. Math.}, 35(3):281--322, 1913.

\bibitem[Mac17]{MR1576566}
P.~A. MacMahon.
\newblock Two applications of general theorems in combinatory analysis: (1) to
  the theory of inversions of permutations; (2) to the ascertainment of the
  numbers of terms in the development of a determinant which has amongst its
  elements an arbitrary number of zeros.
\newblock {\em Proc. London Math. Soc.}, S2-15(1):314, 1917.

\bibitem[{Maz}08]{mazur}
{Mazur, Barry}.
\newblock Bernoulli numbers and the unity of mathematics: Bartlett lecture
  notes, 2008.
\newblock Online.
  \url{http://www.math.harvard.edu/~mazur/papers/slides.Bartlett.pdf}.

\bibitem[MPP18]{1512.08348}
Alejandro~H. Morales, Igor Pak, and Greta Panova.
\newblock Hook formulas for skew shapes {I}. {$q$}-analogues and bijections.
\newblock {\em J. Combin. Theory Ser. A}, 154:350--405, 2018.

\bibitem[MW47]{MR0022058}
H.~B. Mann and D.~R. Whitney.
\newblock On a test of whether one of two random variables is stochastically
  larger than the other.
\newblock {\em Ann. Math. Statistics}, 18:50--60, 1947.

\bibitem[NS11]{Novak-Sniady.2011}
Jonathan Novak and Piotr \'Sniady.
\newblock What is {$\dots$} a free cumulant?
\newblock {\em Notices Amer. Math. Soc.}, 58(2):300--301, 2011.

\bibitem[{OEI}17]{oeis}
{OEIS Foundation Inc.}
\newblock The {O}n-{L}ine {E}ncyclopedia of {I}nteger {S}equences, 2017.
\newblock Online. \url{http://oeis.org}.

\bibitem[RSW04]{Reiner-Stanton-White.CSP}
V.~Reiner, D.~Stanton, and D.~White.
\newblock The cyclic sieving phenomenon.
\newblock {\em J. Combin. Theory Ser. A}, 108(1):17--50, 2004.

\bibitem[Sac97]{MR1453118}
Vladimir~N. Sachkov.
\newblock {\em Probabilistic methods in combinatorial analysis}, volume~56 of
  {\em Encyclopedia of Mathematics and its Applications}.
\newblock Cambridge University Press, Cambridge, 1997.
\newblock Translated from the Russian, Revised by the author.

\bibitem[Sel44]{MR0018287}
Atle Selberg.
\newblock Remarks on a multiple integral.
\newblock {\em Norsk Mat. Tidsskr.}, 26:71--78, 1944.

\bibitem[Sta89]{Stanley.1989}
Richard~P. Stanley.
\newblock Log-concave and unimodal sequences in algebra, combinatorics, and
  geometry.
\newblock In {\em Graph theory and its applications: {E}ast and {W}est
  ({J}inan, 1986)}, volume 576 of {\em Ann. New York Acad. Sci.}, pages
  500--535. New York Acad. Sci., New York, 1989.

\bibitem[Sta99]{ec2}
R.~P. Stanley.
\newblock {\em Enumerative combinatorics. {V}ol. 2}, volume~62 of {\em
  Cambridge Studies in Advanced Mathematics}.
\newblock Cambridge University Press, Cambridge, 1999.

\bibitem[Ste89]{stembridge89}
John~R. Stembridge.
\newblock On the eigenvalues of representations of reflection groups and wreath
  products.
\newblock {\em Pacific J. Math.}, 140(2):353--396, 1989.

\bibitem[SW98]{MR1692145}
John~R. Stembridge and Debra~J. Waugh.
\newblock A {W}eyl group generating function that ought to be better known.
\newblock {\em Indag. Math. (N.S.)}, 9(3):451--457, 1998.

\bibitem[SW10]{Shareshian-Wachs.2010}
John Shareshian and Michelle~L. Wachs.
\newblock Eulerian quasisymmetric functions.
\newblock {\em Adv. Math.}, 225(6):2921--2966, 2010.

\bibitem[Swa18]{s17}
Joshua~P. Swanson.
\newblock On the existence of tableaux with given modular major index.
\newblock {\em Algebraic Combinatorics}, 1(1):3--21, 2018.

\bibitem[TW18]{Theil-Williams-2018}
Marko {Thiel} and Nathan {Williams}.
\newblock {Strange Expectations and the Winnie-the-Pooh Problem}.
\newblock {\em Preprint arXiv:1811.02550}, Nov 2018.

\bibitem[Wik17]{wiki:bernoulli}
Wikipedia.
\newblock Bernoulli number --- {W}ikipedia{,} {T}he {F}ree {E}ncyclopedia,
  2017.
\newblock [Online; accessed 21-July-2017 ].

\bibitem[Zab03]{math/0310301}
Mike Zabrocki.
\newblock A bijective proof of an unusual symmetric group generating function.
\newblock {\em Preprint arXiv:math/0310301}, 2003.

\end{thebibliography}
\bibliographystyle{alpha}

\end{document}